\documentclass[12pt]{article}
\usepackage{graphicx,amsmath,amsthm,amssymb,dsfont, mathrsfs,MnSymbol}
\def\RR{\mathbb R}

\def\ZZ{\mathbb Z}
\def\PP{\mathbb P}

\def\QQ{\mathbb Q}

\def\DD{\mathbb D}

\def\VV{\mathbb V}
\def\UU{\mathbb U}
\def\SS{\mathbb S}

\def\cQ{\mathcal Q}

\def\cH{\mathcal H}

\def\cF{\mathcal F}

\def\fn{\mathfrak n}

\def\cP{\mathcal P}

\def\sD{\mathscr D}
\def\cD{\mathcal D}
\def\dD{\mathds D}

\def\rD{\rm D}
\def\cN{\mathcal N}

\def\fT{\mathfrak T}

\def\fU{\mathfrak U}

\def\bE{{\bf E}}

\def\bx{{\bf x}}
\def\by{{\bf y}}

\def\bm{{\bf m}}

\def\bb{{\bf b}}

\def\br{{\bf r}}

\def\balpha{\boldsymbol{\alpha}}

\def\btheta{\boldsymbol{\theta}}

\def\b1{{\bf 1}}
\def\d1{\mathds{ 1}}
\def\dk{\mathds{ k}}
\def\mod{{\; \rm mod \;}}
\def\ad{{\rm and}}

\def\with{{\rm with}}
\def\where{{\rm where}}
\def\for{{\rm for}}

\def\vol{{\rm vol}}

\def\exp{{\rm exp}}

\def\qed{ \ \vrule width.2cm height.2cm depth0cm\smallskip}

\begin{document}



\title{Temporal Central Limit Theorem for Multidimensional Adding Machine }
\author{Mordechay B. Levin}
\date{}
\maketitle
\begin{abstract}
Let $p_1,...,p_{s+1}$ be distinct primes and let $T_{p_i}$ be the von Niemann - Kakutani adding machine
 $(1 \leq i \leq s)$, $T_{\cP}(\bx) =(T_{p_1}(x_1),..., T_{p_s}(x_s))$. Let $y_i \in (0,1)$ be a $p_{s+1}$-rational
 $(1 \leq i \leq s)$, $\d1_{[0,\by)}$ the indicator function of the box
 $[0,y_1) \times \cdots\times [0,y_s)$. In this paper, we prove the following central limit theorem :
\begin{equation}     \nonumber
  \frac{ \sum_{k=-n}^{n-1}  \d1_{[0,\by)}(T^k_P(\bx)) -2n y_1y_2\dots y_s  }{\cH_N(\bx) \log_2^{s/2} N}
   \; \stackrel{w}{\longrightarrow} \;\cN(0,1),
\end{equation}
when $n$ is sampled uniformly from $\{ 1,...,N\}$, $\cH_N(\bx)  \in [\upsilon_1, \upsilon_2]$
with  some $\upsilon_1, \upsilon_2 >0$,
 for almost all $\bx \in [0,1)^s$.

\end{abstract}
%
%
\section{Introduction }

 Let $\cP_{N} = (\beta_{n})_{n = 0}^{N-1}$ be an $N$-element point set in the
 $s$-dimensional unit cube $[0,1)^s$. The local {\bf discrepancy} function of $\cP_{N}$ is defined as
\begin{equation}\label{In1}
   D(\bx, \cP_{N}  )= \sum\nolimits_{k=0}^{N-1}   \b1_{B_{\bx}}(\beta_{n}) - N x_1 \cdots x_s,
\end{equation}
where $\b1_{B_{\bx}}(\by) =1, \; {\rm if} \;\by  \in B_{\bx}$,
and $   \b1_{B_{\bx}}(\by) =0,$  if $
\by   \notin B_{\bx}$ with $B_{\bx}=[0,x_1) \times \cdots\\ \times [0,x_s) $, $\bx =(x_1,...,x_s)$.
We define the $\emph{L}_q$  discrepancy of $\cP_{N}$ as
\begin{equation} \nonumber 
       D_{s,\infty} (\cP_{N}) =
    \sup_{ 0<x_1, \ldots , x_s \leq 1} \; |
    D(\bx,\cP_{N}) |,\qquad  D_{s,q}(\cP_{N})=\left\| D(\bx,\cP_{N}) \right\|_{s,q}  ,
\end{equation}
\begin{equation*} 
     \left\| f(\bx) \right\|_{s,q}= \Big(  \int_{[0,1]^s}  |f(\bx)|^q d \bx \Big)^{1/q}.
\end{equation*}

 A sequence $(\beta_k)_{k\geq 0}$ is of {\it
low discrepancy} (abbr. l.d.s.) if  $D_{s,\infty}
((\beta_k)_{k=0}^{N-1})=O(\log^s N) $ for $ N \rightarrow \infty $.
An ergodic transformation $T$ of $s-$torus is of
 low discrepancy in the orbit of $\bx$ if  $(T^k(\bx))_{k \geq 0}$ is the l.d.s.

Let  $p_1,...,p_{s}$ be distinct primes
 \begin{equation} \nonumber 
 k=\sum_{j\geq 1}e_{j,i}(k) p_i^{j-1},\;  e_{j,i}(k) \in \{0,1, \ldots
 ,p_i-1\}, \;  {\rm and}    \; \phi_i(k)= \sum_{j\geq 1}e_{j,i}(k) p_i^{-j}.
 \end{equation}
Van der Corput    proved that $ (\phi_1(k))_{k\geq 0}$ is the $1-$dimensional l.d.s.\\
The first example of multidimensional l.d.s. was proposed by Halton :
\begin{equation} \nonumber 
  H_s(k)= (\phi_{1}(k),\ldots ,\phi_{s}(k)), \quad k=0,1,2,... \;.
\end{equation}
For other examples of  l.d.s., see, e.g., [BC],  [Ni].

Let $q \geq 2$ be an integer, and let $x=.x_{1} x_{2} \dots x_{j} \dots =\sum\nolimits_{j=1}^{\infty} x_{j}
 / q^j$  be the
$q$-expansion of $x$.
We define von Neumann-Kakutani's $q-$adic
adding machine:
\begin{equation} \label{In7}
    T_{q}(x) := (x_k+1)/{q}^k + \sum\nolimits_{j \geq k+1}  x_{j} / q^{j}, \qquad \;\;
		T_{q}^n(x)=T_{q}(T_{q}^{n-1}(x)),
\end{equation}
$n=2,3,\dots,T_{q}^0(x)=x,$ where $k= \min \{j \;| \; x_j \ne q-1 \}$.

Let $\cP=(p_1,...,p_{s})$ and let  $T_{\cP}(\bx) =(T_{p_1}(x_1), ..., T_{p_s}(x_s))$.
As is known, the  sequence $(T_{p_i}^k(x))_{n\ge1}$ coincides for $x=0$  with the van
der Corput sequence in base ${p_i}$ $(1 \leq i \leq s)$ (see e.g., [FKP, \S 2.5]).
Hence $T^k_{\cP}(0) =H_s(k)$.
%
It is easy to show that  $(T_{\cP}^k(\bx))_{k \geq 1}$ is the low discrepancy ergodic transformation
for all $\bx$.

 We will consider the probability space $\PP_s =([0,1)^s, B([0,1)^s), \lambda_s)$,
 where $\lambda_s$ is the $s-$dimensional Lebesgue measure on $[0,1)^s$.

Let $T\; : \; [0,1)^s \to  [0,1)^s$ be a map, $f \;: \; [0,1)^s \to \RR $ is a function,
 and $\bx_0 \in [0,1)^s$ is a fixed initial condition.
  We say that the ergodic sums $ S_k(f,\bx_0) =f(\bx_0) + f(T(\bx_0))+ \cdots
+ f(T^{k-1}(\bx_0))$ satisfy a {\it temporal
central limit theorem (TCLT) on the orbit of $\bx_0$}, if there exists
 constants $A_N$ and $B_N$ such that
\begin{equation} \nonumber 
  (S_k -A_N)/B_N   \stackrel{w}{\rightarrow} \cN(0,1),
\end{equation}
when $k$ is sampled uniformly from $\{ 1,...,N\}$.

The ergodic sum $ S_k(f,\bx)$ is said to satisfy an almost sure TCLT on $\PP_s $,
 if for $\lambda_s$-a.e. $\bx$,
 $ S_k(f,\bx)$  satisfy the TCLT.

The first example of the TCLT was discovered by Beck [Be1],[Be2]. Beck considered the transformation $x \to x +
 \sqrt{2} \mod 1$ and the indicator function
of an interval $[0,y)$.
In [DS1], Dolgopyat and Sarig considered a more general notion of the temporal distributional limit
 theorem (see also [ADS], [BrUl], [DS2], [DS3] ).


In [Le2], we proved an $s+1$-parametric spatio-temporal CLT (see definition in [DS1, p.12])
for an $s-$dimensional adding machine (Halton's sequence) with indicator functions.
Similar results for another low discrepancy ergodic transformations was obtained in [Le1] (see also [LeMe]).

Let  $y_i= \sum\nolimits_{j \geq 1} y_{i,j} / p_i^j$  be the
$p_i$-expansion of $y_i$  $(i=1,...,s)$, $\kappa_1, \kappa_2 >0$ and let
\begin{equation}  \label{In10} 
 \liminf_{N \to \infty} \min_{1 \leq i \leq s} \#\{ j \in [1,N]  \; | \;
   1 \leq y_{i,j}\; \; \ad \;\; \{y_i p_i^{j} \} \leq 1 -\kappa_1  \}/N = \kappa_2 >0.
\end{equation}
For example, \eqref{In10} is true for $q$-rational reels $y_i>0$, where $(q,p_i)=1$ $(1 \leq i \leq s)$.
It is easy to see that  \eqref{In10} is also true for a.e. $\by  =(y_1,...,y_s) \in [0,1)^s$. Let
\begin{multline} \label{CLT-1}%
 \sD_{\by,\bx}(L)= D(\by, (T_{\cP}^k(\bx))_{k =- L}^{L-1} ), \;\;\;\;
  \dot{\cH}_N(\bx)=
 \int_{[0,1]}  \sD_{\by,\bx}([\theta N])    d \theta, \\
         \ddot{\cH}_N(\bx)=
 \Big(\int_{[0,1]}  \sD^2_{\by,\bx}([\theta N])    d \theta \Big)^{1/2}  .
\end{multline}
In this paper, we prove the existence of an almost sure TCLT for the discrepancy function in the
 multidimensional case : \\ \\
{\bf Theorem.} {\it Let $s \geq 1$, $\theta$ be a uniformly distributed random variable in $[0,1]$,
$\by$ satisfy \eqref{In10}. Then
\begin{equation}  \nonumber 
    \sD_{\by,\bx}([\theta N]) / \ddot{\cH}_N(\bx) \;
	\stackrel{w}{\longrightarrow} \;\cN(0,1), \qquad \qquad \dot{\cH}_N(\bx)/ \ddot{\cH}_N(\bx) \;
	\stackrel{N \to \infty}{\longrightarrow} \; 0,
\end{equation}
where $ (\log_2 N)^{-s/2} \ddot{\cH}_N(\bx)  \in [\pi^{-1} p_0^{-4-s/2} 2^{-s/2}
\kappa_1^{s} \kappa_2^{s/2},
 7 p_0^{1+s/2}]$ for $N \geq N_0(\bx)$
with some $N_0(\bx)$ $\for \;a.s.\;  \bx \in [0,1)^s$, $p_0 =p_1 p_2 \cdots p_s$.}
\\

Now we describe the structure of the paper. In Lemma 1, we get a simple
estimate of Fourier series of the discrepancy function $\sD_{\by,\bx}(L)$.
In Lemma 2 and Lemma 4, we minimise
the number of terms in the expression of  $\sD_{\by,\bx}(L)$. The main tool of
the proof is the $S$-unit theorem and the theorem on linear forms in logarithm (see \S 2.2).
In \S 2.3 and 2.4, we compute the upper and the lower bounds of the variance of
 $\sD_{\by,\bx}(L)$.

In order to prove the Theorem, we use the martingale CLT (see Theorem D).  We construct a martingale
 approximation to  $\sD_{\by,\bx}(L)$ in \S 2.6.
In order to apply Theorem D, in Lemma 12, we compute the sum of four moment of parts of
 $\sD_{\by,\bx}(L)$ and in Lemma 13, we find Levi's conditional variance of  $\sD_{\by,\bx}(L)$.

\section{Proof of the Theorem}
{\bf 2.1. \;\; Beginning of the proof of the Theorem.}\\

 We will use notation $A \ll B$  equal to $A = O(B)$.
Let
\begin{equation} \nonumber
   \Delta(\fT) =   \begin{cases}
    1,  & \; {\rm if}  \;  \fT  \;{\rm is \;true},\\
    0, &{\rm otherwise},
  \end{cases} \qquad \qquad  \delta_M(a) =   \begin{cases}
    1,  & \; {\rm if}  \;  a \equiv 0  \mod M,\\
    0, &{\rm otherwise}.
  \end{cases}
\end{equation}

Let $[y]$ be the integer part of $y$,
\begin{equation} \label{Beg-1a}
   I_M=[-[(M-1)/2],[M/2]]\cap \ZZ, \qquad\qquad  I_M^{*}=I_M \setminus \{ 0 \}.
\end{equation}
 Note that  the integers of the interval $I_M$ are a complete set of residues $\mod M$, $M \geq 1$.
By [Ko, Lemma 2, p. 2], we have
\begin{equation} \label{Beg-1}
   \delta_M(a) =   \frac{1}{M} \sum_{k\in I_M} e\Big(\frac{ak}{M}\Big),  \qquad \where \quad e(x) =\exp(2\pi i x).
\end{equation}
By [Ko, Lemma 1, p. 1], we get
\begin{equation} \label{Beg-2a}
  \Big|\frac{1}{M} \sum_{k=0}^{M-1} e\big(k \alpha) \Big|
  \leq \min\big(1, \frac{1}{2M\llangle \alpha \rrangle} \big) ,\quad \;\; \where \;
	\llangle \alpha \rrangle = \min( \{\alpha\}, 1 - \{\alpha\}).
\end{equation}
Let $\bar{m} = \max(1, |m|) \leq M/2 $.     From [Ko, p. 2], we obtain for $R \leq M$ :
\begin{equation} \label{Beg-2}
  \Big|\frac{1}{M} \sum_{k=0}^{R-1} e\big(\frac{m k}{M}\big) \Big| \leq \min\Big(1, \Big|\frac{e(m R/M) -1}{M
	(e(m /M)-1)} \Big|\Big) \leq \frac{|\sin(\pi m R/M)|}{\bar{m}}  \leq \frac{1}{\bar{m}}.
\end{equation}
Let $x_i =0.x_{i,1}x_{i,2}...=\sum_{j \geq 1}  x_{i,j} p_{i}^{-j}$, with $x_{i,j} \in \{0,1,...,p_{i}-1  \}$,
 $i=1,...,s$.
We define the truncation
\begin{equation}  \nonumber
        [x_i]_r =\sum_{1 \leq j \leq r} x_{i,j} p_{i}^{-j} \quad \with \quad r \geq 1.
\end{equation}
If $\bx = (x_1, . . . , x_s)  \in [0, 1)^s$, then the truncation $[\bx]_{\br}$ is defined coordinatewise,
 that is, $[\bx]_{\br}=
( [x_1]_{r_1}, . . . , [x_s]_{r_s})$, where $\br =(r_1,...,r_s)$.

Let
\begin{equation} \label{Beg-2b}
   V_{i,r}(x_i)=
  \sum_{1 \leq j \leq r}  x_{i,j} p_i^{j-1} \;  \in [0,p_i^{r_i}.
\end{equation}
By (\ref{In7}), we have
\begin{equation}\nonumber
  V_{i,r}(T_{p_i}(x_i)) \equiv V_{i,r}(x_i) +1  \; {\rm mod} \; p_{i}^{r}.
\end{equation}
Hence
\begin{equation}\nonumber
  V_{i,r}(T^k_{p_i}(x_i)) \equiv V_{i,r}(x_i) +k \; {\rm mod} \; p_{i}^{r}, \quad k=0,1,...
\end{equation}
and
\begin{equation} \label{Beg-3a}
    [T^k_{p_i}(x_i)]_{r}=[y_i]_{r} \;
			\Leftrightarrow  \; k \equiv  V_{i,r}(y_i)- V_{i,r}(x_i)      \; ({\rm mod} \; p_{i}^{r}),
          \quad 1 \leq i \leq s .
\end{equation}
 Let  $p_0=p_1 p_2\cdots p_s$, $ P_{\br} =p_1^{r_1}p_2^{r_2} \dots p_s^{r_s}$
  and let $  M_{i,\br}$ be the unique
integer satisfying the two conditions
\begin{equation} \label{Beg-3} 
   M_{i,\br} \equiv
	\big( P_{\br}/p_i^{r_i} \big)^{-1} \mod p_i^{r_i}, \; M_{i,\br} \in [0,p_i^{r_i}), \; 1 \leq i \leq s.
\end{equation}
We define $ V_{\br,\bx} \in [0,P_{\br})$ as follows
\begin{equation}\label{Beg-6} 
  V_{\br,\bx} \equiv \sum_{i=1}^s  M_{i,\br}
   P_{\br}  p_i^{-r_i}  V_{i,r_i}(x_i) \equiv  \sum_{i=1}^s  M_{i,\br}
   P_{\br}  p_i^{-r_i}    \sum_{1 \leq j \leq r_i}  x_{i,j} p_i^{j-1} \;
    \mod \; P_{\br}.
\end{equation}
Applying (\ref{Beg-3a}) and the Chinese Remainder Theorem, we get
\begin{equation} \label{Beg-10} 
 [T^k_{\cP}(\bx)]_{\br} = [\by]_{\br}  \; \Longleftrightarrow
 \; k \equiv   V_{\br,\by}  - V_{\br,\bx}  \; ({\rm mod} \; P_{\br}).
\end{equation}
Let $n =[\log_2 N]+1$, $\theta \in [0,1)$, $L= [\theta N]$.
From [Ni, p. 29, 30] and \eqref{CLT-1}, we get
\begin{equation}  \label{Beg-16}
  \cD(L) :=  \sD_{[\by]_n,\bx}(L) =  \sD_{\by,\bx}(L)  )
   + \epsilon s,  \;\;\; |\epsilon| \leq 1.
\end{equation}
Let
\begin{equation}  \label{Beg-18}
 V_{\br,\by,\bb} \equiv \sum_{i=1}^s  M_{i,\br}
   P_{\br}  p_i^{-r_i}
  \Big( \sum_{1 \leq j < r_i}  y_{i,j} p_i^{j-1}  + b_i p_i^{r_i-1} \Big) \mod  P_{\br} , \quad
	 V_{\br,\by,\bb} \in [0,P_{\br}).
\end{equation}
Similarly to [Ni, p. 37-39], we obtain from \eqref{In1}, \eqref{CLT-1}, \eqref{Beg-10} and \eqref{Beg-18}
\begin{equation}  \label{Beg-20}
 \cD(L)=   \sum_{k=-L}^{L-1} \sum_{r_1,...,r_s =1}^n \sum_{b_1=0}^{y_{1,r_1}-1} \cdots  \sum_{b_s=0}^{y_{s,r_s} -1} \delta_{P_{\br}} (k -  V_{\br,\by,\bb} + V_{\br,\bx})- 2L [y_1]_n \cdots [y_s]_n .
\end{equation}
Let
\begin{equation}  \label{Beg-18a}
  \dD_{\br,L} =  \sum_{b_1=0}^{y_{1,r_1}-1} \cdots  \sum_{b_s=0}^{y_{s,r_s} -1}  \sum_{k=-L}^{L-1} (\delta_{P_{\br}} (k  - V_{\br,\by,\bb} + V_{\br,\bx})-1/P_{\br} ).
\end{equation}
It is easy to see that $ | \dD_{\br,L}| \leq  y_{1,r_1}\cdots y_{s,r_s}  < p_0$ and
\begin{equation}  \label{Beg-22}
     \cD(L)  =\sum_{r_1,...,r_s =1}^n   \dD_{\br,L}.
\end{equation}
Let  $W_0= 50 p_0 s^2 \log_2 n$,
\begin{equation}   \label{Beg-26}
  U= \{ (r_1,...,r_2) \in [1,n]^s \;  | \;  \max_i(r_i) >W_0\}, \; \;\;
             \tilde{\cD}(L)=  \sum_{\br \in U}    \dD_{\br,L}.
\end{equation}
From \eqref{Beg-20} - \eqref{Beg-26}, we derive
\begin{equation}  \nonumber 
| \cD(L) -  \tilde{\cD}(L) | \leq  p_0 W_0^s = p_0^{s+1} (50 s^2)^s \log_2^{s} n  \ll \log_2^{s} n.
\end{equation}
Using  \eqref{Beg-16}, we obtain
\begin{equation}  \label{Beg-28}
  \sD_{\by,\bx}(L)   =   \tilde{\cD}(L)    + O(\log_2^{s} n).
\end{equation}\\
{\bf Lemma 1.} {\it With the notations as above, we have}
\begin{multline}  \label{Le1-1}
\dD_{\br,L} =
    \sum_{m \in I_{P_{\br}}^{*}}   \varphi_{ \br,L,m} \; \psi_{ \br}(m,\by)\;
e\Big( \frac{m}{P_{\br}} (V_{\br,\bx} - V_{\br,\by}) \Big), \;\;  \varphi_{ \br,L,m}
  =\frac{2 i  \sin(2\pi m L /P_{\br})}{P_{\br}(e(m/P_{\br})-1)}, \qquad\qquad  \\
    |\varphi_{ \br,L,m}| \leq \frac{|\sin(2\pi m L /P_{\br}) |}{\bar{m}} \leq \frac{1}{\bar{m}}, \;
  \psi_{ \br}(m,\by)  = \prod_{i=1}^s \dot{\psi}(i,\{-m M_{i,\br}/p_i \}p_i, y_{i, r_i}) , \\
\dot{\psi}(i,m', y_{i, r_i}) = \sum_{0 \leq b  < y_{i, r_i}} e(m'(b-y_{i, r_i} )/p_i),\;\;
 \;\;
				\dot{\psi}(i,m', y_{i, r_i}) = \frac{1- e(-m' y_{i, r_i}/p_i)}{e(m'/p_i)-1} \\
	 \for\;\; \{m'/p_i\} \neq 0  , \quad  \dot{\psi}(i,0,y_{i, r_i})=y_{i, r_i},\quad  \dot{\psi}(i,m', 0)=0,
 \quad
  |\dot{\psi}(i,m', y_{i, r_i})| \leq p_i .
\end{multline} \\
{\bf Proof.}
Using \eqref{Beg-1a}, \eqref{Beg-1} and \eqref{Beg-18a}, we obtain
\begin{multline*}
 \dD_{\br,L}  =    \sum_{b_1=0}^{y_{1,r_1} -1} \cdots  \sum_{b_s=0}^{y_{s,r_s} -1} \sum_{k=-L}^{L-1} \frac{1}{P_{\br}}
 \sum_{m \in I_{P_{\br}}^{*}} e\Big( \frac{m}{P_{\br}} (k - V_{\br,\by,\bb}+ V_{\br,\bx})\Big) \\
=  \sum_{b_1=0}^{y_{1,r_1} -1} \cdots  \sum_{b_s=0}^{y_{s,r_s} -1}  \sum_{m \in I_{P_{\br}}^{*}}
\frac{e(m L/P_{\br}) - e(-m L /P_{\br} )}{P_{\br}(e(m/P_{\br})-1)}
\; e\Big( \frac{m}{P_{\br}} (V_{\br,\bx} -V_{\br,\by,\bb}) \Big) \\
= \sum_{m \in I_{P_{\br}}^{*}} e\Big( \frac{m}{P_{\br}}
 (V_{\br,\bx} -V_{\br,\by})  \Big) \varphi_{ \br,L,m}
\sum_{b_1=0}^{y_{1,r_1} -1} \cdots  \sum_{b_s=0}^{y_{s,r_s} -1}
e\Big(  \frac{m}{P_{\br}}( V_{\br,\by} -V_{\br,\by,\bb} ) \Big).
\end{multline*}
According to \eqref{Beg-6} and  \eqref{Beg-18}, we get
\begin{equation*}  
 V_{\br,\by,\bb} \equiv V_{\br,\by} + \sum_{i=1}^s  M_{i,\br}
   P_{\br}  (b_i - y_{i,r_i})/p_i  \mod  P_{\br}.
\end{equation*}
%
Therefore
\begin{equation*}
 \dD_{\br,L} = \sum_{m \in I_{P_{\br}}^{*}} e\Big( \frac{m}{P_{\br}}
 (V_{\br,\bx} -V_{\br,\by})  \Big) \varphi_{ \br,L,m} \;\varpi_0,
\end{equation*}
with
\begin{equation*}
\varpi_0= \sum_{b_1=0}^{y_{1,r_1} -1} \cdots  \sum_{b_s=0}^{y_{s,r_s} -1}
e\Big(  -m\sum_{i=1}^s  M_{i,\br}(b_i-y_{i,r_i})/p_i \Big).
\end{equation*}
It is easy to see that
\begin{multline*}
\varpi_0= \prod_{i=1}^s  e(  m M_{i,\br}  y_{i,r_i}/p_i ) \sum_{b_i=0}^{y_{i,r_i} -1}
e(  -m M_{i,\br} b_i /p_i )  \\
 = \prod_{i=1}^s  \frac{1- e(m M_{i,\br} y_{i, r_i}/p_i)}{e(-m M_{i,\br}/p_i)-1}
=\prod_{i=1}^s \dot{\psi}(i,\{-m M_{i,\br}/p_i \}p_i, y_{i, r_i}) =\psi_{ \br}(m,\by).
\end{multline*}
Hence
\begin{equation*}
   \dD_{\br,L} =    \sum_{m \in I_{P_{\br}}^{*}}   \varphi_{ \br,L,m} \; \psi_{ \br}(m,\by) \;
e\Big( \frac{m}{P_{\br}}  (V_{\br,\bx} -V_{\br,\by}) \Big).
\end{equation*}
Bearing in mind that $\sin(\pi x/2) \geq x$ for $x \in [0,1]$, we obtain from
\eqref{Beg-2} and \eqref{Le1-1} that
\begin{align*}
& |\varphi_{ \br,L,m}  | = \Big| \frac{2\sin(2\pi m L /P_{\br})}{P_{\br}(e(m/P_{\br})-1)} \Big|
 = \frac{|\sin(2\pi m L /P_{\br}) |}{|P_{\br} \sin(\pi m  /P_{\br})}
 \leq \frac{|\sin(2\pi m L /P_{\br}) |}{2 \bar{m}}
     \leq \frac{1}{2\bar{m}},      \\
&     \dot{\psi}(i,m', y_{i, r_i})=
 \sum_{0 \leq b < y_{i, r_i}} e(m(b-y_{i, r_i})/p_i),
 \;\; |\dot{\psi}(i,m', y_{i, r_i})| \leq p_i, \quad |\psi_{ \br}(m',\by)| \leq p_0. \nonumber
\end{align*}
If $ y_{i, r_i}= 0$ then $\dot{\psi}(i,m', y_{i, r_i})=0 $. Let $ y_{i, r_i} \neq 0$.
 If $ \{m'/p_i\}=0$, then $ \dot{\psi}(i,m', y_{i, r_i})=y_{i, r_i} \geq 1$.
 Hence Lemma 1 is proved. \qed \\ \\
{\bf Lemma 2.} {\it Let $\balpha =(\alpha_1,...,\alpha_s) \in \ZZ^s$,
\begin{multline}  \label{Le2-1}
 \dot{\dD}_{\br,L} = \sum_{ \substack{\alpha_i \in [0, 10s \log_2 n]\\ 1 \leq i \leq s}}
    \sum_{\substack{ m P_{\balpha} \in I_{P_{\br}}^{*} \\ (m,p_0)=1 }}   \varphi_{ \br,L,m P_{\balpha}} \;
     \psi_{ \br}(m P_{\balpha},\by)\;
e\Big( \frac{m P_{\balpha}}{P_{\br}}  (V_{\br,\bx} -V_{\br,\by}) \Big), \\
                \dot{\cD}(L)=  \sum_{\br \in U}     \dot{\dD}_{\br,L}.
\end{multline}
Then}
\begin{equation} \nonumber
    \dD_{\br,L} =  \dot{\dD}_{\br,L} + O(n^{-9 s}), \quad   \dD_{\br,L} \ll n \qquad \ad \quad
     \tilde{\cD}(L) =  \dot{\cD}(L) + O(n^{- 8s}),
\end{equation} \\
{\bf Proof.} Let $m = m^{'}P_{\balpha}$, with $(m^{'}, p_0)=1$. Using Lemma 1, we get
\begin{equation} \nonumber 
  \dD_{\br,L} = \sum_{ \substack{\alpha_i \geq 0\\ 1 \leq i \leq s}} \;
    \sum_{\substack{  m^{'} P_{\balpha} \in I_{P_{\br}}^{*} \\ ( m^{'},p_0)=1 }}
      \varphi_{ \br,L, m^{'} P_{\balpha}} \; \psi_{ \br}( m^{'} P_{\balpha},\by)\;
e\Big( \frac{ m^{'}P_{\balpha}}{P_{\br}}  (V_{\br,\bx} -V_{\br,\by}) \Big).
\end{equation}
Hence
\begin{multline}  \nonumber 
 | \dD_{\br,L} -\dot{\dD}_{\br,L}| \leq
 \sum_{ \substack{\max_i \alpha_i > 10s \log_2 n\\
  \alpha_i \geq 0, 1 \leq i \leq s}}
    \sum_{\substack{ m^{'} P_{\balpha} \in I_{P_{\br}}^{*} \\ (m^{'},p_0)=1 }}
    \big|  \varphi_{ \br,L,m^{'} P_{\balpha}} \; \psi_{ \br}(m^{'} P_{\balpha},\by) \big| \\
 \leq p_0 \sum_{ \substack{\max_i \alpha_i > 10s \log_2 n\\
  \alpha_i \geq 0, 1 \leq i \leq s}}
    \sum_{\substack{ m^{'} P_{\balpha} \in I_{P_{\br}}^{*} \\ (m^{'},p_0)=1 }}
    ( \bar{m}^{'} P_{\balpha}   )^{-1}
    \ll n \sum_{ \substack{\max_i \alpha_i > 10s \log_2 n\\
  \alpha_i \geq 0, 1 \leq i \leq s}} P_{-\balpha} \ll n^{-9s}.
\end{multline}
It is easy to see that $ \dD_{\br,L} \ll  \sum_{|m| \leq p_0^n} 1/\bar{m}   \ll n$.
Bearing in mind that $\# U \leq n^s$,
we obtain from \eqref{Beg-26} and \eqref{Le2-1} the assertion of Lemma 2. \qed \\ \\

We consider  expectations in  probability spaces $\PP_1$ and $\PP_s$ :
\begin{equation}  \label{Le2-5a}
     \bE_{\btheta}f_1 = \bE_{\btheta}(f_1)= \int_{[0,1)} f_1(\theta)d \theta , \quad          \bE_{x}f_2=  \bE_{x}(f_2)= \int_{[0,1)^s} f_2(\bx)d \bx.
\end{equation}
It is easy to see that
\begin{equation} \label{Le2-5}
     \bE_{\btheta}f([\theta N]) =  \frac{1}{N}  \sum_{k=0}^{N-1} f( k).
\end{equation}
Bearing in mind that $\sin(z)=(\exp(iz)-\exp(-iz))/2i$, $\cos(z) =\sin(z+\pi/2) $, $2\sin^2(z) = 1-\cos(2z)= 1-\sin(2z+\pi/2)$ and  $2\sin(z_1) \sin(z_2) = \cos(z_1-z_2) - \cos(z_1+z_2) =
\sin(z_1-z_2 +\pi/2) - \sin(z_1+z_2+\pi/2) $, we obtain from  \eqref{Beg-2a}, \eqref{Le2-5}
\begin{align} \label{Le2-8}
&    |\bE_{\btheta}\sin (2\pi [\theta N] \alpha + \beta)| = N^{-1} \Big| \sum_{k=0}^{N-1} \sin (2\pi k \alpha + \beta))   \Big|
   \leq \min\Big(1, \frac{1}{2N\llangle \alpha \rrangle} \Big), \nonumber  \\
& 2|\bE_{\btheta}\sin^2 (2\pi [\theta N] \alpha) -1/2|   \leq \min\Big(1, \frac{1}{2N\llangle 2\alpha \rrangle} \Big),\\
&   2 |\bE_{\btheta}\sin (2\pi  [\theta N] \alpha) \sin (2\pi [\theta N] \alpha)|
     \leq \min\Big(1, \frac{1}{2N\llangle \alpha -\beta \rrangle} \Big) +
   \min\Big(1, \frac{1}{2N\llangle \alpha +\beta\rrangle} \Big). \nonumber
\end{align} \\
{\bf Lemma 3.}  {\it Let $\mu \geq 2$ be integer, $ \br_0 =(r_{0,1},...,r_{0,s}), \; \;
r_{0,i} =\max_{1 \leq j \leq \mu}(r_{j,i})$,
\begin{equation}  \nonumber 
  \omega_1 =  \sum_{j=1}^{\mu} m_j (V_{\br_j,\bx} - V_{\br_j,\by}  )  /P_{\br_j} .
\end{equation}
Then }
\begin{equation}   \nonumber 
     \bE_{x}e(  \omega_1)
 =    \delta_{P_{\br_0} }\big( \sum_{j=1}^{\mu} m_j P_{\br_0-\br_j} \big), \quad \ad \quad
 \bE_{x}e(  \omega_1) =\Delta\big( \sum_{j=1}^{\mu} m_j P_{-\br_j} =0 \big)
\end{equation}
 for  $ |\sum_{j=1}^{\mu} m_j P_{\br_0-\br_j}| < P_{\br_0} $. \\  \\
{\bf Proof.}
By \eqref{Beg-1}, \eqref{Beg-2b},  \eqref{Beg-3}, \eqref{Beg-6} and \eqref{Le2-5a}, we obtain
\begin{equation}   \nonumber 
     V_{\br,\bx} /P_{\br}  =
    \sum_{i=1}^s  M_{i,\br}
     \sum_{ \nu=1}^{ r_{i}}   \frac{x_{i,\nu} p_i^{\nu-1} }{ p_i^{r_i} }   \equiv
      \sum_{i=1}^s  \frac{M_{i,\br}}{ p_i^{r_i} }
     \sum_{ \nu=1}^{ r_{0,i}}  x_{i,\nu} p_i^{\nu-1} \equiv
       \sum_{i=1}^s  V_{i,r_{0,i}}(x_i)  \frac{M_{i,\br}}{ p_i^{r_i} }\; \mod 1
\end{equation}
and
\begin{equation}  \nonumber 
\bE_{x}e( m V_{\br,\bx} /P_{\br})    =\bE_{x}e( m \sum_{i=1}^s  M_{i,\br}
     \sum_{ \nu=1}^{ r_{i}}   x_{i,\nu} p_i^{\nu-1} / p_i^{r_i}) =
  \prod_{i=1}^s \delta_{p_i^{r_i}} (m  M_{i,\br}) =
\delta_{P_{\br}} (m).
\end{equation}
Hence
\begin{multline}   \nonumber 
 \sum_{j=1}^{\mu} m_j   (V_{\br_j,\bx}-V_{\br_j,\by} ) /P_{\br_j} \equiv
    \sum_{j=1}^{\mu} m_j    \sum_{i=1}^s  M_{i,\br_j}  ( V_{i,r_{j,i}}(x_i)- V_{i,r_{j,i}}(y_i) )
    /p_i^{r_{j,i}} \\
\equiv
    \sum_{j=1}^{\mu} m_j    \sum_{i=1}^s  M_{i,\br_j}  ( V_{i,r_{0,i}}(x_i)- V_{i,r_{0,i}}(y_i) )
    /p_i^{r_{j,i}} \\
\equiv
  \sum_{i=1}^s ( V_{i,r_{0,i}}(x_i)- V_{i,r_{0,i}}(y_i) ) /p_i^{r_{0,i}} \sum_{j=1}^{\mu} m_j     M_{i,\br_j}
   p_i^{r_{0,i} -r_{j,i}}    \; \mod 1.
\end{multline}
Let $ \omega_2 =e\big( -  \sum_{1 \leq i  \leq s} V_{i,r_{0,i}}(y_i)  /p_i^{r_{0,i}} \sum_{1 \leq j \leq \mu} m_j
   M_{i,\br_j}
   p_i^{r_{0,i} -r_{j,i}}     \big)$. We have
\begin{multline}   \nonumber 
 \bE_{x}e(  \omega_1)
 =  \omega_2 \prod_{i=1}^s \delta_{p_i^{r_{0,i}} } \big( \sum_{j=1}^{\mu} m_j  M_{i,\br_j}
   p_i^{r_{0,i} -r_{j,i}}   \big)=  \prod_{i=1}^s \delta_{p_i^{r_{0,i}} } \big( \sum_{j=1}^{\mu} m_j  M_{i,\br_j}
   p_i^{r_{0,i} -r_{j,i}}   \big)  \\
   = \prod_{i=1}^s \delta_{p_i^{r_{0,i}} } \big( \prod_{1 \leq \nu \leq s, \nu \neq i}
	p_{\nu}^{ r_{0,\nu}}
   \sum_{j=1}^{\mu} m_j   M_{i,\br_j}
   p_i^{r_{0,i} -r_{j,i}}   \big).
\end{multline}
From  \eqref{Beg-3},  we get
\begin{multline}   \nonumber 
   M_{i,\br_j} \equiv \prod_{1 \leq {\nu} \leq s, {\nu} \neq i}
	p_{\nu}^{-r_{j,\nu}}  \mod p_i^{r_{j,i}} \\
 \ad  \qquad
               p_i^{r_{0,i} -r_{j,i}}\prod_{1 \leq {\nu} \leq s, {\nu} \neq i}
	p_{\nu}^{ r_{0,\nu}} M_{i,\br_j}
\equiv \prod_{1 \leq {\nu} \leq s}
 p_{\nu}^{r_{0,{\nu}} -r_{j,\nu}}
  \mod p_i^{r_{0,i}}.
\end{multline}
Therefore
\begin{multline}    \nonumber 
 \bE_{x}e(  \omega_1)   =   \prod_{i=1}^s \delta_{p_i^{r_{0,i}} } \big(
\sum_{j=1}^{\mu} m_j    p_{i}^{r_{0,i} -r_{j,i}}    \prod_{1 \leq {\nu} \leq s}
 p_{\nu}^{r_{0,{\nu}} -r_{j,\nu}}  M_{\nu,\br_j}    \big) \\
 =   \prod_{i=1}^s \delta_{p_i^{r_{0,i}} } \big(
\sum_{j=1}^{\mu} m_j \prod_{1 \leq {\nu} \leq s}
 p_{\nu}^{r_{0,{\nu}} -r_{j,\nu}}    \big)=     \delta_{P_{\br_{0}} } \big(
\sum_{j=1}^{\mu} m_j \prod_{1 \leq {\nu} \leq s}
 p_{\nu}^{r_{0,{\nu}} -r_{j,\nu}}    \big) \\
 =     \delta_{P_{\br_{0}} } \big(
\sum_{j=1}^{\mu} m_j \prod_{1 \leq {\nu} \leq s}
 p_{\nu}^{r_{0,{\nu}} -r_{j,\nu}}    \big)
  =    \delta_{P_{\br_0} }\big( \sum_{j=1}^{\mu} m_j P_{\br_0-\br_j} \big).
\end{multline}
Now consider the case  $ |\sum_{j=1}^{\mu} m_j P_{\br_0-\br_j}| < P_{\br_0}$.
It is easy to see that, if
$   \sum_{j=1}^{\mu} m_j P_{\br_0-\br_j} \equiv 0 \mod P_{\br_0}  $, then
$   \sum_{j=1}^{\mu} m_j P_{\br_0-\br_j}=0$.

Hence Lemma 3 is proved. \qed \\ \\
{\bf Lemma 4.}  {\it Let
\begin{multline}
 \ddot{\dD}_{\br,L} = \sum_{ \substack{\alpha_i \in [0, 10s \log_2 n]\\ 1 \leq i \leq s}} \;
    \sum_{\substack{ 0 < |m| \leq n^{10s} \\ (m,p_0)=1 }}   \varphi_{ \br,L,m P_{\balpha}} \;
    \psi_{ \br}(m P_{\balpha},\by)\;
e\Big( \frac{m P_{\balpha}}{P_{\br}}  (V_{\br,\bx} -V_{\br,\by}) \Big), \\
             \ddot{\cD}(L)=  \sum_{\br \in U}     \ddot{\dD}_{\br,L}.  \label{Le3-1}
\end{multline}
Then there exists  $N_1(\bx)$ such that
\begin{align}
&    \bE_{\btheta} \big(( \sD_{\by,\bx}([\theta N])  ) -
       \ddot{\cD}([\theta N]))^2\big)  \ll \log_2^{2s} n,   \label{Le3-2}\\
&   0.5 \bE_{\btheta} \ddot{\cD}^2([\theta N]) - \log_2^{2s} n     \leq  \bE_{\btheta}  \sD^2_{\by,\bx}([\theta N])
           \leq 2\bE_{\btheta} \ddot{\cD}^2([\theta N]) + \log_2^{2s} n, \nonumber
\end{align}
with  $ N > N_1(\bx)$ for  a.s. $ \bx \in [0,1)^s$.} \\  \\
{\bf Proof.}
Let $a= \sD_{\by,\bx}([\theta N])  )$, $b=\dot{\cD}([\theta N]) $,
$c=\ddot{\cD}([\theta N]) $.
 By  \eqref{Beg-28} and Lemma 2,   we get
\begin{equation*}
      \sD_{\by,\bx}([\theta N])  ) = \tilde{\cD}([\theta N]) +O(\log_2^{s} n),
       \;\;
      \tilde{\cD}([\theta N]) =\dot{\cD}([\theta N])+O(n^{-8s}).
\end{equation*}
Therefore $  a-b =O(\log_2^{s} n) $.
%
It is easy to see that
\begin{multline*}  %
     \bE_{\btheta} (a-c)^2 \leq 2  \bE_{\btheta} (a-b)^2 + 2 \bE_{\btheta} (b-c)^2 \\
ad \quad  0.5 \bE_{\btheta} c^2 - 2  \bE_{\btheta} (a-b)^2 - 2 \bE_{\btheta} (b-c)^2
    \leq  0.5 \bE_{\btheta} c^2 -  \bE_{\btheta} (a-c)^2  \leq  \bE_{\btheta} a^2   .
\end{multline*}
Hence
\begin{multline*}  
    \bE_{\btheta} ( \sD_{\by,\bx}([\theta N])   -
       \ddot{\cD}([\theta N]))^2) \\
     \qquad\qquad\qquad\qquad\qquad   \leq    2
  \bE_{\btheta} (\dot{\cD}([\theta N]) -
       \ddot{\cD}([\theta N]))^2 + \log_2^{2s+1} n  \\
\ad \;  0.5 \bE_{\btheta} \ddot{\cD}^2([\theta N]) -  \bE_{\btheta} (\dot{\cD}([\theta N]) -
       \ddot{\cD}([\theta N]))^2 - \log_2^{2s+1} n  \leq
          \bE_{\btheta}  \sD^2_{\by,\bx}([\theta N]) ) .
\end{multline*}
%
Bearing in mind that
\begin{equation}  \nonumber 
  \bE_{\btheta} a^2 \leq 2  \bE_{\btheta} (a-c)^2 +  2\bE_{\btheta} c^2
 \leq
 2 \bE_{\btheta} c^2  +  4\bE_{\btheta} (a-b)^2 + 4 \bE_{\btheta} (b-c)^2,
\end{equation}
we have
\begin{equation}  \nonumber 
          \bE_{\btheta}  \sD^2_{\by,\bx}([\theta N]) ) \leq
  2\bE_{\btheta} \ddot{\cD}^2([\theta N]) +  4\bE_{\btheta} (\dot{\cD}([\theta N]) -
       \ddot{\cD}([\theta N]))^2 + \log_2^{2s+1} n
\end{equation}
for $N \geq N_2$ with some $N_2>1$.\\
Therefore, in order to prove Lemma 4 it is enough to verify that
\begin{equation}  \nonumber 
          \bE_{\btheta}  |\dot{\cD}([\theta N]) - \ddot{\cD}([\theta N]) |^2 \ll 1
          \qquad   \for \; a.s. \; \bx \in [0,1)^s.
\end{equation}
Applying the Borel-Cantelli lemma, we get that \eqref{Le3-2} results from the following inequality
\begin{equation}  \nonumber 
        {\rm Prob}\Big( \bx \in [0,1)^s \; | \;
         \bE_{\btheta} \big(\dot{\cD}([\theta N]) - \ddot{\cD}([\theta N]) \big)^2 > 1 \Big) \ll  1/n^2.
\end{equation}
By Chebyshev's inequality, it is enough to prove that
\begin{equation}   \label{Le3-6}
   \bE_{x}\bE_{\btheta} \big|\dot{\cD}([\theta N]) - \ddot{\cD}([\theta N]) \big|^2   \ll  1/n^2.
\end{equation}
From  \eqref{Le2-1} and \eqref{Le3-1}, we have
\begin{multline}  \label{Le3-7}
     \bE_{x}  \bE_{\btheta}   \big| \dot{\cD}([\theta N]) - \ddot{\cD}([\theta N]) \big|^2 =
     \sum_{\br_1,\br_2 \in U} Z_{\br_1,\br_2}, \;\; \with  \;\;   Z_{\br_1,\br_2}:=
\sum_{ \substack{\alpha_{j,i} \in [0, 10s \log_2 n]\\ 1 \leq i \leq s, 1 \leq j \leq 2}}\\
   \times \sum_{\substack{ |m_j|> n^{10s} \\ (m_j,p_0)=1, j=1,2 }}
    \bE_{\btheta} \Big( \varphi_{ \br_1,[\theta N],m_1 P_{\balpha_1}}
      \overline{ \varphi_{ \br_2,[\theta N],m_2 P_{\balpha_2}}} \Big)
       \psi_{ \br_1}(m_1 P_{\balpha_1},\by) \overline{\psi_{ \br_1}(m_2 P_{\balpha_2},\by)} \\
 \times
 \bE_{x}  e\Big( \frac{m_1 P_{\balpha_1}}{P_{\br_1}}  (V_{\br_1,\bx} -V_{\br_1,\by}) -
  \frac{m_2 P_{\balpha_2}}{P_{\br_2}}  (V_{\br_2,\bx} -V_{\br_2,\by})   \Big) .
\end{multline}
Using  Lemma 1 and Lemma 3, we obtain
\begin{multline} \label{Le3-8}
   Z_{\br_1,\br_2} =
 \sum_{ \substack{\alpha_{j,i} \in [0, 10s \log_2 n]\\ 1 \leq i \leq s, 1 \leq j \leq 2}}
    \sum_{\substack{ |m_j|> n^{10s} \\ (m_j,p_0)=1, j=1,2 }}
    \bE_{\btheta} \Big( \varphi_{ \br_1,[\theta N],m_1 P_{\balpha_1}}  \overline{ \varphi_{ \br_2,[\theta N],m_2 P_{\balpha_2}}} \Big) \\
 \times
    \; \psi_{ \br_1}(m_1 P_{\balpha_1},\by) \overline{\psi_{ \br_1}(m_2 P_{\balpha_2},\by)}
 \delta_{P_{\br_0} }( P_{\br_0-\br_1+\balpha_1}m_1 - P_{\br_0-\br_2+\balpha_2}m_2 ) \\
  \ll
\sum_{ \substack{\alpha_{j,i} \in [0, 10s \log_2 n]\\ 1 \leq i \leq s, 1 \leq j \leq 2}}
  \sum_{\substack{ |m_j|> n^{10s} \\ (m_j,p_0)=1, j=1,2 }} (\bar{m}_1 \bar{m}_2 P_{\balpha_1}  P_{\balpha_2}  )^{-1} \\
  \times   \delta_{P_{\br_0} }( P_{\br_0-\br_1+\balpha_1}m_1 - P_{\br_0-\br_2+\balpha_2}m_2 ),
\end{multline}
with $ \br_0 =(r_{0,1},...,r_{0,s}), \; \; r_{0,i} =\max(r_{1,i},r_{2,i}), \; i=1,...,s $.

By \eqref{Le3-7} and \eqref{Le3-8}, we get  that \eqref{Le3-6} go after the following inequality
\begin{multline} \label{Le3-9}
 \varpi_1:=
\sum_{ \substack{\alpha_{j,i} \in [0, 10s \log_2 n]\\ 1 \leq i \leq s, 1 \leq j \leq 2}}
  \sum_{\substack{ |m_j|> n^{10s} \\ (m_j,p_0)=1, j=1,2 }}
   (\bar{m}_1 \bar{m}_2 P_{\balpha_1}  P_{\balpha_2}  )^{-1} \\
  \times   \delta_{P_{\br_0} }( P_{\br_0-\br_1+\balpha_1}m_1 - P_{\br_0-\br_2+\balpha_2}m_2 )
    \ll n^{-2s-2} .
\end{multline}
From  \eqref{Beg-26}, we have  that there exists $i_0 \in [1,s]$
such that $r_{1,i_0 } \geq W_0 =50s^2 p_0 \log_2 n$.  Let $\rho = \max_{j=1,2} r_{j,i_0}$,
 $\rho_j = r_{j,i_0} -\alpha_{j,i_0}$,
 $\rho_{j_1} = \max_{j=1,2} \rho_j$ with $j_1 \in \{1,2\}$.
Bearing in mind that $m_1,m_2$ satisfy the congruence
$P_{\br_0-\br_1+\balpha_1}m_1 - P_{\br_0-\br_2+\balpha_2}m_2 \equiv 0 \mod  P_{\br_0}$ and that
$\max_{j,i} \alpha_{j,i} \leq 10s \log_2 n$,  we get $ \rho_{j_1 } \geq  40s^2  \log_2 n$ and
\begin{equation*}
      m_{j_1} p_{i_0}^{\rho -\rho_{j_1} }    \equiv    m_{j_2} p_{i_0}^{\rho -\rho_{j_2} } A_0
      \mod p_{i_0}^{\rho}, \quad \with \quad j_2  \in \{1,2\} \setminus \{j_1\}
\end{equation*}
for some integer $A_0$ with $  (A_0,p_{i_0})=1$.
Hence
\begin{equation*}
      m_{j_1}   \equiv    m_{j_2} p_{i_0}^{\rho_{j_1} -\rho_{j_2} } A_0
      \mod p_{i_0}^{\rho_{j_1}}, \quad \quad  \rho_{j_1 } \geq  40 s^2  \log_2 n.
\end{equation*}
Therefore
\begin{multline*}
 \varpi_1 \leq  \sum_{n^{10s} <  |m_{j_2}| \leq 2^n  } \frac{1}{\bar{m}_{j_2}}
  \sum_{n^{10s} <  |m_{j_1}| \leq 2^n  } \frac{ \delta_{p_{i_0}^{[40 s^2  \log_2 n]}}
   (  m_{j_1}   -   m_{j_2} p_{i_0}^{\rho_{j_1} -\rho_{j_2} } A_0 )  }{\bar{m}_{j_1}}\\
 \ll \sum_{n^{10s} <  |m_{j_2}| \leq 2^n  } \frac{1}{\bar{m}_{j_2}}
  \sum_{0 \leq  \ell \leq 2^n  } \frac{ 1 }{n^{10s}
    + \ell p_{i_0}^{ 40 s^2  \log_2 n }}
    \ll n^{-9s}.
\end{multline*}
Hence \eqref{Le3-9}, \eqref{Le3-6}, \eqref{Le3-2} and Lemma 4 are proved. \qed \\ \\

{\bf 2.2. Diophantine inequalities.}

 We consider the following simple variant of the \texttt{ S-unit theorem} (see  [ESS,  Theorem~1.1, p. 808]):
Let $\beta_1,...,\beta_d \in \QQ$, $\beta_i \neq 0, \; i=1,...,d$. We consider the equation
\begin{equation}\label{d-1}
      \beta_1 P_{\br_1} +... +\beta_d P_{\br_d} =1, \quad P_{\br_j}=\prod_{i=1}^s p_i^{r_{j,i}},\;\; r_{j,i} \in \ZZ.
\end{equation}
  A solution $(\br_1,...,\br_d)$ of (\ref{d-1})
 is called \texttt{ non-degenerate}
if  $ \sum_{i \in J}  \beta_i P_{\br_i}  \neq 0$ for every nonempty subset $J$ of $\{1,...,d\}$.  \\ \\
{\bf Theorem  A.} {\it The number  $A(\beta_1,...,\beta_d)$ of non-degenerate solutions
  of equation (\ref{d-1})
satisfies the estimate }
\begin{equation}  \nonumber 
     A(\beta_1,...,\beta_d) \leq \exp((6d)^{3d}s).
\end{equation}
\\
 {\bf Corollary 1.} {\it  Let $\balpha_j \in \ZZ^s$,  $|m_j| >0$, $1 \leq j \leq 2d$, $d=2,4$.  Then}
\begin{equation}  \nonumber 
 \sum_{\substack{ \br_j \in\UU  \\ 1 \leq j \leq 2 d}}
 \Delta \Big(  \sum_{1 \leq j \leq 2 d}  m_j P_{- \br_j+\balpha_j} =0\Big)
   \ll \#\UU^{d}.
\end{equation}
%
\texttt{   Linear forms in logarithms.}
Write $\Lambda$ for the linear form in logarithms,
\begin{equation}\nonumber
    \Lambda =b_1 \ln \alpha_1 + ...+ b_{\tilde{k}} \ln \alpha_{\tilde{k}}     ,
\end{equation}
where $b_1,...,b_{\tilde{k}}$ are integers, $|b_i| \leq B\; (i=1,...,{\tilde{k}})$, $B \geq e$. We shall assume  that $\alpha_1,...,\alpha_{\tilde{k}} $
 are non-zero algebraic numbers
with  heights at most $A_1,...,A_{\tilde{k}}$ (all $\geq e$) respectively.  \\ \\
 {\bf Theorem B.} [BW, Theorem 2.15, p. 42] {\it  If $\Lambda \neq 0$, then
\begin{equation}\nonumber 
   | \Lambda |  > \exp( -(16{\tilde{k}}d)^{2({\tilde{k}}+2)} \ln A_1 ... \ln A_{\tilde{k}}  \ln B) ,
\end{equation}
where $d$ denote the degrees of} $\QQ(\alpha_1,...,\alpha_{\tilde{k}} )$.
\\ \\
 {\bf Corollary 2.} Let $|r_i| \leq 2n$, $0 < |m_j| \leq n^{10s}$, $|m_1/m_2|P_{\br} > 1$. Then
\begin{equation} \nonumber 
    |m_1/m_2P_{\br} -1|  > \exp( -C_1 \ln^2 n) , \quad \with \quad C_1 =20s(16(s+1))^{2(s+3)} \ln^s p_0.
\end{equation} \\
{\bf Proof.} Taking into account that $exp(x) -1 \geq x$ for $x \geq 0$, we get
\begin{equation} \nonumber
   | |m_1/m_2|P_{\br} -1| \geq    \ln(|m_1/m_2|P_{\br} )=  r_1 \ln(p_1) +...+ r_s \ln(p_s) + \ln(|m_1/m_2 |).
\end{equation}
Using Theorem B with $\tilde{k}=s+1$, $b_i=r_i$, $\alpha_i =A_i=p_i$, $(i=1,...,s)$,
$b_{s+1} =1$, $\alpha_{s+1} =1$, $A_{s+1} = n^{10s}$, $B=2n$, $d=1$, we get the assertion of the Corollary 1. \qed \\ \\

{ \bf 2.3. Upper bound of the variance of $ \ddot{\cD}_1([\theta N])$ }.\\

By \eqref{Beg-26},  $W_0= 50 s^2 p_0 \log_2 n $ and $U= \{ (r_1,...,r_s) \in [1,n]^s \;  | \;  \max_i r_i >W_0\}$.
Let
\begin{align}  \label{d-10}
&  W_1=[n/ W_2]-2, \quad  W_2= [\log_2^{20s} n],\quad \;\; W_3= [\log_2^{10s} n], \quad \;\;
  A_k = n- (k+1)W_2,  \nonumber\\
&  B_k =A_k + W_2 - 2W_3, \qquad  h(\br)  =r_1\log_2(p_1)+ \cdots + r_s \log_2(p_s),\nonumber \\
&    U_1 =  \bigcup_{1 \leq k \leq W_1 } \UU_{k},\quad  \UU_{k}=
 \{\br \in U \; | \;  h(\br) \in [A_k,B_k) \}, \nonumber \\
&    U_2 =  \bigcup_{1 \leq k \leq W_1 } \UU_{k}^{'},\quad  \UU_{k}^{'}=
 \{\br \in U \; | \;  h(\br) \in [B_{k+1},A_k) \}, \nonumber \\
&  U_3 = \{\br \in U \; | \;   h(\br)  \in [B_{1}, n+W_2]  \}, \qquad
  U_4 = \{\br \in U \; | \;   h(\br)  >n+W_2 \}, \nonumber \\
&   U_5 = \{\br \in U \; | \;   h(\br)  <B_{W_1+1} \},\qquad \ddot{\cD}_j(L)=  \sum_{\br \in U_j}
   \ddot{\dD}_{\br,L}, \quad
{\DD}_k(L)=  \sum_{\br \in \UU_{k}}    \ddot{\dD}_{\br,L}.
\end{align}
According to Lemma 1 and Lemma 4, we have
\begin{equation}  \label{d-11}
     {\DD}_k(L) \ll n^{s+1}, \qquad L \in [1,N].
\end{equation}
We see that $  U = \bigcup_{1 \leq j \leq 5 } U_j, \;\; U_i \cap U_j =\emptyset, \;
 i \neq j, \; 1 \leq i,j \leq 5 $.\\
By  \eqref{Le3-1} and \eqref{d-10},  we get that
\begin{equation}  \label{d-12}
    \ddot{\cD}(L)=  \sum_{1 \leq j \leq 5 }   \ddot{\cD}_j(L) \qquad \ad \qquad
     \ddot{\cD}_1(L)=  \sum_{ 1 \leq k \leq W_1}     {\DD}_k(L).
\end{equation}  \\
 {\bf Theorem C.} [Wi, p.59, Theorem 2, ref. 3] {\it Let $\dot{\Gamma} \subset \RR^s$ be a lattice,
  $\det \dot{\Gamma} =1$, $\cQ \subset \RR^s$ a compact convex body
 and $\dot{r}> \sqrt{s}/2$ the radius of
its greatest sphere in the interior. Then}
\begin{equation}\nonumber
       \vol( \cQ) \big(1- 0.5 \sqrt{s}/\dot{r}  \big)  \leq  \# \dot{\Gamma}
        \cap \cQ  \leq  \vol( \cQ )\big(1+ 0.5 \sqrt{s}/\dot{r}   \big) .
\end{equation} \\
Let $\fU \subset \RR^s$, $\fU \dot{-} \balpha = \{\br \in \ZZ^s \; | \; \br + \balpha \in \fU\} $,
\begin{multline}\label{d-30}
\fU_{\balpha_1,\balpha_2}=   \{\br \in \ZZ^s \;|\; \br+\balpha_1 \in \fU, \br +\balpha_2  \in \fU \}
  =  \{\br \in \fU \dot{-} \balpha_1 \} \cap \{\br \in \fU \dot{-} \balpha_2  \} , \\
\hat{\fU}_{\balpha_1,\balpha_2} = \fU \cap \fU_{\balpha_1,\balpha_2}, \quad
      \dot{\fU}_{\balpha_1,\balpha_2} = \fU \setminus \hat{\fU}_{\balpha_1,\balpha_2}, \quad
\ddot{\fU}_{\balpha_1,\balpha_2} = \fU_{\balpha_1,\balpha_2} \setminus \hat{\fU}_{\balpha_1,\balpha_2}.
\end{multline}
It is easy to see that
\begin{equation} \label{d-32}
 \fU_{\balpha_1,\balpha_2}= ( \fU \setminus  \dot{\fU}_{\balpha_1,\balpha_2}) \bigcup
 \ddot{\fU}_{\balpha_1,\balpha_2}, \quad   \dot{\fU}_{\balpha_1,\balpha_2} \in \fU, \quad \ad \quad
     \ddot{\fU}_{\balpha_1,\balpha_2} \notin \fU. \\
\end{equation} \\
{\bf Lemma 5.}   {\it Let $ U_6 =[0,2W_2]^s \cup U_2 $. With the notations as above}
\begin{multline} \nonumber 
     U_{1_{\balpha_1,\balpha_2}}= ( U_1 \setminus  \dot{U}_{1_{\balpha_1,\balpha_2}}) \bigcup
 \ddot{U}_{1_{\balpha_1,\balpha_2}}, \quad
  \dot{U}_{1_{\balpha_1,\balpha_2}} \in U_1, \quad \ad \quad \ddot{U}_{1_{\balpha_1,\balpha_2}} \notin U_1,\\
\# \UU_k \ll n^{s-1} \log_2^{20s} n, \qquad  \#U_2 \leq n^s \log_2^{-10s} n,\qquad
  \; \#U_6 \leq n^s \log_2^{-10s} n,\\
\#U_3 \leq n^{s-1} \log_2^{40s} n,
 \#\dot{U}_{1_{\balpha_1,\balpha_2}} \;+ \;\#\ddot{U}_{1_{\balpha_1,\balpha_2}}
\ll n^{s} \log_2^{-19s} n,    \;\;\;
    \;\;\;
  \# U_5 \ll \log_2^{21s^2} n,  \\
       \#U_1 \leq n^s , \qquad    U_1 \supset  [1, n/(p_0 s)]^s \setminus U_6.
\end{multline} \\
{\bf Proof.}
The first three assertions follow from \eqref{d-32}. \\
\textsf{ Consider} $\UU_k=  \{\br \in U \; | \;  h(\br) \in [A_k,B_k) \} $  ($ A_k = n- (k+1)W_2$,
$B_k =A_k + W_2 - 2W_3$). We have that if $g \leq  h(\br)  =r_1\log_2(p_1)+ \cdots + r_s \log_2(p_s)\leq f$,
then $g/\log_2 p_0 \leq \max_i r_i \leq f$, where $p_0 =p_1 p_2 \cdots p_s$.
It is easy to see that $\UU_k $ is included in a rotated and shifted parallelepiped
 $[0,sn]^{s-1} \times [0,2W_2-4W_3]$ :
\begin{equation*} 
   \UU_k  \subset M_k [0,sn]^{s-1} \times [0,2W_2-4W_3] + z_k, \quad \for \quad 1 \leq k \leq W_1=[n/W_2],
\end{equation*}
$   W_2 =[\log_2^{20s} n]$,  with some $M_k \in  SO(s)$ and $z_k \in \RR^s$,
where  $ SO(s)$ is the special orthogonal group.
Applying Theorem C, we obtain
\begin{equation*} 
   \# \UU_k \ll  n^{s-1} W_2  \quad \ad \quad
   \# U_3 = \sum_{k \in [B_1, n+W_2)} \# \UU_k  \ll n^{s-1} W_2^2 \ll n^{s-1} \log_2^{40s} n.
\end{equation*}
\textsf{ Consider} $\UU_k^{'} $.  It is easy to see that there exist
$M_k^{'} \in  SO(s)$ and $z_k^{'} \in \RR^s$ such that
\begin{equation*} 
   \UU_k^{'} \subset M_k^{'} [0,sn]^{s-1} \times [0,2W_3] + z_k^{'}.
\end{equation*}
Applying Theorem C, we get
\begin{multline*} 
   \# \UU_k^{'} \ll  n^{s-1} W_3 \ll n^{s-1} \log_2^{10s} n, \quad \qquad \#U_2 \leq
    W_1 \max_k \# \UU_k^{'} \ll n^s \log_2^{-10s} n, \\
    \ad  \qquad \#U_6 \ll W_2^s +   \# U_2 \ll \log_2^{20s^2} + \# U_2 \ll n^s \log_2^{-10s} n.
\end{multline*}
\textsf{ Consider} $ U_{1_{\balpha_1,\balpha_2}}$.  It is easy to verify that
\begin{equation*} 
   U_{1_{\balpha_1,\balpha_2}} = \bigcup_{k=1}^{W_1} \UU_{k,\balpha_1,\balpha_2}, \quad
   \dot{U}_{1_{\balpha_1,\balpha_2}} = \bigcup_{k=1}^{W_1} \dot{\UU}_{k,\balpha_1,\balpha_2},
    \;\;\; \ad \;\;\; \ddot{U}_{1_{\balpha_1,\balpha_2}} = \bigcup_{k=1}^{W_1}  \ddot{\UU}_{k,\balpha_1,\balpha_2}.
\end{equation*}
We have that there exist
$M_{k,j} \in  SO(s)$, $z_{k_j}\in \RR^s$  and $\mu_0 \leq 2^{2s}$ such that
\begin{equation*} 
  \dot{\UU}_{k,\balpha_1,\balpha_2} \subset  \bigcup_{1 \leq j
   \leq \mu_0} M_{k,j} [0,sn]^{s-1} \times [0,20 s \log_2 n] + z_{k,j}.
\end{equation*}
Applying Theorem C, we obtain
\begin{equation*} 
   \# \dot{\UU}_{k,\balpha_1,\balpha_2} \ll  n^{s-1} \log_2 n , \quad \ad \quad
    \#\dot{U}_{1_{\balpha_1,\balpha_2}} \ll  n^{s} \log_2^{-19s} n .
\end{equation*}
Similarly, we get that  $  \# \dot{U}_{2_{\balpha_1,\balpha_2}} \ll  n^{s} \log_2^{-19s} n $.\\
Using \eqref{d-10}, we get $ \# U_5 \ll W_2^s \ll \log_2^{20s^2} n $, $ \# U_1 \leq n^s $ and
 $U_1 \supset  [1, n/(p_0 s)]^s \setminus U_6$.
%
%
%
%
Hence  Lemma 5 is proved. \qed \\ \\
Let
\begin{multline}  \label{d-16}
 \chi_{1,\bm,\br} = \Delta(  m_1 P_{\balpha_1 -\br_1} =  -m_2 P_{\balpha_2 -\br_2} ),  \qquad
   \chi_{2,\bm,\br} = \Delta(  m_1 P_{\balpha_1 -\br_1}  = m_2 P_{\balpha_2 -\br_2}),\\
   \chi_{3,\bm,\br} =1-  \chi_{1,\bm,\br}-  \chi_{2,\bm,\br},
\end{multline}
and let
\begin{multline} \label{d-18}
\rD_{\nu, \mu} = \sum_{\br_1, \br_2 \in U_{\nu}}
    \sum_{ \substack{\alpha_{j,i} \in [0, 10s \log_2 n]\\ 1 \leq i \leq s, j=1,2}}
    \sum_{\substack{  0 < |m_j| \leq n^{10s} \\ (m_j,p_0)=1,j=1,2 }}
      \varphi_{ \br_1,[\theta N],m_1 P_{\balpha_1}}
      \varphi_{ \br_2,[\theta N],m_2 P_{\balpha_2}}
     \\
\times
    \psi_{ \br_1}(m_1 P_{\balpha_1},\by) \psi_{ \br_2}(m_2 P_{\balpha_2},\by)\\
    \times
e\Big( \frac{m_1 P_{\balpha_1}}{P_{\br_1}}  (V_{\br_1,\bx} -V_{\br_1,\by})   +  \frac{m_2 P_{\balpha_2}}{P_{\br_2}}
 (V_{\br_2,\bx} -V_{\br_2,\by}) \Big) \; \chi_{\mu,\bm,\br}.
\end{multline}
Using  \eqref{Le3-1}, we have
\begin{equation}  \label{d-20}
\ddot{\cD}_{\nu}^2([\theta N])) =   \rD_{\nu,1} +\rD_{\nu,2} +   \rD_{\nu,3}.
\end{equation} \\
{\bf Lemma 6.}   {\it With the notations as above}
\begin{align*} 
 &  \ddot{\cD}_4([\theta N])  \ll 1, \qquad \quad  \ddot{\cD}_5([\theta N])  \ll \log_2^{21s^2}  n,
         \qquad \quad  \bE_{\btheta}  \ddot{\cD}_1([\theta N]) \ll 1/n,\\
&   \bE_{\btheta} \ddot{\cD}_1^2([\theta N]) \leq 32 p_0^{s+2} n^{s} ,
    \quad  \bE_{\btheta} \ddot{\cD}_2^2([\theta N])  \ll n^{s} \log^{-10s} n ,
    \quad
  \bE_{\btheta} \rD_{1,3} \ll 1. \nonumber
\end{align*} \\
{\bf Proof.}
\textsf{ Consider}  $ \ddot{\cD}_4([\theta N])$. Using \eqref{Le3-1}, \eqref{d-10} and  Lemma 1, we get
\begin{equation}  \nonumber
 |\ddot{\cD}_4([\theta N])|  \leq   \sum_{ \br \in U_4}  \sum_{ \substack{\alpha_i \in [0, 10s \log_2 n]\\ 1 \leq i \leq s}}
    \sum_{\substack{  0 < |m| \leq n^{10s} \\ (m,p_0)=1 }}
       \frac{|\sin(2\pi m [\theta N]/P_{\br-\balpha})|}{\bar{m}P_{\balpha}} \;
          |\psi_{ \br}(m P_{\balpha},\by)|.
\end{equation}
By \eqref{d-10}, we have
 $P_{\br} > 2^{n+W_2} \geq N 2^{\log_2^{20s} n -1}$.
Bearing in mind that $|\sin(x)| \leq |x|$, we obtain
\begin{multline} \label{d-118} 
 |\ddot{\cD}_4([\theta N])|     \leq         \sum_{ \br \in U_4}   \sum_{ \substack{\alpha_i \in [0, 10s \log_2 n]\\ 1 \leq i \leq s}}
    \sum_{\substack{  0 < |m| \leq n^{10s} \\ (m,p_0)=1 }}
   \frac{2\pi N}{P_{\br}}p_0  \\
    \times \Delta(P_{\br} >  N 2^{\log_2^{20s} n -1})
          \ll
      n^s\; \log_2^s n \;n^{10s} \; n^{-300s} \ll 1.
\end{multline}
%
%
\textsf{ Consider}  $ \ddot{\cD}_5([\theta N])$. From \eqref{Le3-1}, \eqref{d-10} and  Lemma 1, we derive
\begin{multline*}  
 \ddot{\cD}_5([\theta N])     \leq        p_0 \sum_{ \br \in U_5}   \sum_{ \substack{\alpha_i \in [0, 10s \log_2 n]\\ 1 \leq i \leq s}}
    \sum_{\substack{  0 < |m| \leq n^{10s} \\ (m,p_0)=1 }}
    \frac{1}{\bar{m}P_{\balpha}} \ll  \sum_{ \br \in U_5}  \log_2 n \\
          \ll  \sum_{ \max_i r_i \leq 2s W_2}  \log_2 n \ll  W_2^s    \log_2 n  \ll   \log_2^{20s^2+1} n.
\end{multline*}
%
\textsf{ Consider}  $\bE_{\btheta}  \ddot{\cD}_1([\theta N])$. Using  \eqref{Le2-8}, \eqref{Le3-1}, \eqref{d-10} and  Lemma 1, we obtain
\begin{multline}  \nonumber
 |\bE_{\btheta} \ddot{\cD}_1([\theta N])|  \leq   \sum_{ \br \in U_1}
  \sum_{ \substack{\alpha_i \in [0, 10s \log_2 n]\\ 1 \leq i \leq s}} \;
    \sum_{\substack{  0 < |m| \leq n^{10s} \\ (m,p_0)=1 }}
       \frac{|\bE_{\btheta}\sin(\frac{2\pi m [\theta N]}{P_{\br-\balpha}})|}{\bar{m}P_{\balpha}} \;
          |\psi_{ \br}(m P_{\balpha},\by)|  \\
\leq p_0    \sum_{ \br \in U_1}  \sum_{ \substack{\alpha_i \in [0, 10s \log_2 n]\\ 1 \leq i \leq s}}
    \sum_{\substack{  0 < |m| \leq n^{10s} \\ (m,p_0)=1 }}
       \frac{1}{\bar{m}P_{\balpha}} \;
          \big(1,  \frac{1}{2N\llangle m  P_{-\br+\balpha}    \rrangle} \big).
\end{multline}
Taking into account that $ |m| \leq n^{10s}  $ and
$ 2^{\log^{20s}_2 n} \leq P_{\br}  \leq 2^{n -\log^{10s}_2 n} $ for $\br \in U_1$, we get
\begin{equation}  \nonumber
 |\bE_{\btheta} \ddot{\cD}_1([\theta N])|
\ll    \sum_{ \br \in U_1}  \sum_{ \substack{\alpha_i \in [0, 10s \log_2 n]\\ 1 \leq i \leq s}}
    \sum_{\substack{  0 < |m| \leq n^{10s} \\ (m,p_0)=1 }}
       \frac{1}{\bar{m}^2 P_{\balpha}} \;
          \frac{  P_{\br - \balpha}}{N}  \ll n^{s}  2^{- \log^{10s}_2 n}  \ll 1/n.
\end{equation}
\textsf{ Consider}  $ \rD_{\nu,\mu}$  $(\nu, \mu=1,2)$.  Suppose that $\chi_{\mu,\bm,\br} \geq 1 $.
Bearing in mind that $(m_j,p_0)=1,\; j=1,2$, we
have  from  \eqref{d-16} that $ |m_1|=|m_2| $ and $\br_2 =\br_1 +\balpha_2 -\balpha_1 $.
Applying \eqref{d-18}, we obtain
\begin{multline*} 
 | \rD_{\nu, 0}| \leq  \sum_{\br_1, \br_2 \in U_{\nu}}
  \sum_{ \substack{\alpha_{j,i} \in [0, 10s \log_2 n]\\ 1 \leq i \leq s, j=1,2}}
    \sum_{\substack{  0 < |m_j| \leq n^{10s} \\ (m_j,p_0)=1,j=1,2 }}
 \frac{p_0^2}{\bar{m}_1 \bar{m}_2 P_{\balpha_1 + \balpha_2}}
\Delta( |m_1|=|m_2|) \\
  \times   \Delta(\br_2 =\br_1 +\balpha_2 -\balpha_1)
  \leq \# U_{\nu}     \sum_{ \substack{\alpha_{j,i} \in [0, 10s \log_2 n]\\ 1 \leq i \leq s, j=1,2}}
    \sum_{\substack{  0 < |m_1| \leq n^{10s} \\ (m_1,p_0)=1 }}
 \frac{2 p_0^2}{\bar{m}_1^2 P_{\balpha_1 + \balpha_2}} \leq 32 p_0^{s+2} \# U_{\nu}.
\end{multline*}
From Lemma 5,  we get
\begin{equation}   \label{Le5-31}
  \rD_{1, \mu} \leq 32 p_0^{s+2} n^s, \qquad \qquad \rD_{2, \mu} \ll n^s /\log_2^{20s} n, \quad \mu=1,2.
\end{equation}
\textsf{ Consider}  $\rD_{\nu,3}$ $(\nu=1,2)$. Applying \eqref{d-18} and Lemma 1, we get
\begin{multline*}
  | \bE_{\btheta} \rD_{\nu,3}| \ll  \sum_{\br_1, \br_2 \in U_{\nu}}
   \sum_{ \substack{\alpha_{j,i} \in [0, 10s \log_2 n]\\ 1 \leq i \leq s, j=1,2}}
    \sum_{\substack{  0 < |m_j| \leq n^{10s} \\ (m_j,p_0)=1,j=1,2 }}
    \Big| \bE_{\btheta} \Big( \varphi_{ \br_1,[\theta N],m_1 P_{\balpha_1}} \nonumber \\
 \times    \varphi_{ \br_2,[\theta N],m_2 P_{\balpha_2}}
   \Big) \Big|
  \chi_{3,\bm,\br}
 \ll  \sum_{\br_1, \br_2 \in U_{\nu}}     \sum_{ \substack{\alpha_{j,i} \in [0, 10s \log_2 n]\\
 1 \leq i \leq s, j=1,2}}
    \sum_{\substack{  0 < |m_j| \leq n^{10s} \\ (m_j,p_0)=1,j=1,2 }} 1 \nonumber \\
  \times | P_{\br_1}(e(m_1/P_{\br_1 -\balpha_1})-1) |^{-1}
     | P_{\br_2}(e(m_2/P_{\br_2 -\balpha_2})-1) |^{-1}   \\ 
 \times  \Big| \bE_{\btheta}   \big(  \sin( 2\pi m_1 [\theta N] /P_{\br_1 -\balpha_1})
  \sin( 2\pi m_2 [\theta N] /P_{\br_2 -\balpha_2})     \big) \Big|
  \chi_{3,\bm,\br}  .  \nonumber
\end{multline*}
Using \eqref{Le2-8}, we derive
\begin{multline}  \label{Le4-7}
 | \bE_{\btheta} \rD_{\nu,3}|
 \ll      \sum_{\br_1, \br_2 \in U_{\nu}}
 \sum_{ \substack{\alpha_{j,i} \in [0, 10s \log_2 n]\\ 1 \leq i \leq s, j=1,2}}
    \sum_{\substack{  0 < |m_j| \leq n^{10s} \\ (m_j,p_0)=1,j=1,2 }}
  \frac{ 1}{\bar{m}_1 \bar{m}_2 P_{\balpha_1 +\balpha_2}} \\
      \times \sum_{\ell=-1,1}\min\big(1,
    \frac{1}{2N\llangle m_1 P_{-\br_1+\balpha_1} + \ell m_2 P_{-\br_2+\balpha_2}   \rrangle} \big)
    \chi_{3,\bm,\br}  .
\end{multline}
By \eqref{d-10}, we get $W_2= [\log_2^{20s} n]$, $\br_{\nu} \in U_{\nu}$,
 $h_{\br_{j}}  \leq n-W_2$ $(\nu,j =1,2) $,
\begin{equation}  \nonumber
  N/P_{\br_j} \geq 2^{W_2} \geq 2^{\log^{20s} n -2} \quad \ad \quad m_1/ P_{\br_1-\balpha_1} \neq {\ell}m_2/ P_{\br_2-\balpha_2}
\end{equation}
for $\ell=-1,1$, $\nu =1,2$ and  $\chi_{3,\bm,\br}  > 0$.
 From \eqref{Beg-26}, we obtain  $  0 <|m_j/P_{\br_j - \alpha_j}| <1/10$ for $\br_{\nu} \in U_{\nu}$.
Using Corollary 2, we have
\begin{multline*}  
 N\llangle m_1 P_{-\br_1+\balpha_1} + m_2 P_{-\br_2+\balpha_2}      \rrangle =
 N | m_1/ P_{\br_1-\balpha_1} + (-1)^{\ell}m_2/ P_{\br_2-\balpha_2} | \\
  \geq \frac{N}{2} \Big( P_{-\br_2+\balpha_2}|m_2| \;
 \big|\frac{|m_1|}{|m_2|}P_{\br_2 -\br_1 +\balpha_1 -\balpha_2} -1\big| +
  P_{-\br_1+\balpha_1}|m_1| \;
 \big|\frac{|m_2|}{|m_1|}P_{\br_1 -\br_2 +\balpha_2 -\balpha_1} -1\big|
 \Big) \\
  \geq N/2 \min\big( P_{-\br_2+\balpha_2}|m_2|, P_{-\br_1 +\balpha_1}|m_1| \big)
  \exp( -C_1 \ln^2 n)   \\
  \geq
 2^{W_2-1}  \exp( -C_1 \ln^2 n)  \geq 2^{\log_2^{20s} n -4}  \exp( -C_1 \ln^2 n)     \gg n^{100s}.
\end{multline*}
We get from  \eqref{Le4-7} that
\begin{equation}  \label{Le4-9}
    | \bE_{\btheta} \rD_{\nu,3}| \ll  \sum_{\br_1, \br_2 \in U_{\nu}}  \;
     \sum_{ \substack{\alpha_{j,i} \in [0, 10s \log_2 n]\\ 1 \leq i \leq s, j=1,2}}
    \sum_{\substack{  0 < |m_j| \leq n^{10s} \\ (m_j,p_0)=1,j=1,2 }}
    \frac{1}{\bar{m}_1 \bar{m}_2} n^{-100s} \ll 1, \quad \nu=1,2.
\end{equation}
By  \eqref{d-20} $\ddot{\cD}_{\nu}^2([\theta N])) =   \rD_{\nu,1} +   \rD_{\nu,2}+   \rD_{\nu,3}$ $(\nu=1,2)$.
Taking into account
\eqref{Le5-31} and \eqref{Le4-9}, we obtain that $ \bE_{\btheta} \ddot{\cD}_1^2([\theta N]) \leq 65 p_0^{s+2}  n^{s}$,
    $\;\;\;\bE_{\btheta} \ddot{\cD}_2^2([\theta N])  \ll n^{s} \log^{-10s} n $,
$ \;\;\;\bE_{\btheta} \rD_{\nu,3} \ll 1 $.
Hence Lemma 6 is proved. \qed \\ \\
%
%
{\bf Lemma 7.} {\it With the notations as above}
\begin{equation*}  
  \bE_{\btheta} \ddot{\cD}_3^2([\theta N]) \ll n^{s -1/8}  \qquad \for  \;\; a.s. \;\; \bx.
\end{equation*} \\
{\bf Proof.}
By \eqref{Le3-1} and \eqref{d-10},  we get
\begin{align}
& \bE_{x}  \bE_{\btheta}^2\ddot{\cD}_3^2([\theta N])   = \sum_{ \substack{ \br_j \in U_3 \\ j \in [1,4] }}
      \sum_{ \substack{\alpha_{j,i} \in [0, 10s \log_2 n]\\ 1 \leq i \leq s,j \in [1,4] }}
    \sum_{\substack{  0 < |m_j| \leq n^{10s} \\ (m_j,p_0)=1,j \in [1,4] }}   \vartheta  \label{Le5-1}   \\
&  \times     \prod_{1 \leq j_1 \leq 2} \bE_{\btheta} \Big(  \prod_{1 \leq j_2 \leq 2}
   \varphi_{ \br_{j_2 +2(j_1-1)},[\theta N],m_{j_2 +2(j_1-1)}  P_{\balpha_{j_2 +2(j_1-1)}}} \Big)
     \prod_{1 \leq j \leq 4} \psi_{ \br_j}(m_j P_{\balpha_j},\by)
   \; \; \nonumber \\
&   \with \quad      \vartheta = \bE_{x}   e\Big(
\sum_{1 \leq j \leq 4}  \frac{m_jP_{\balpha_j }}{P_{\br_j}}
 (V_{\br_j,\bx} -V_{\br_j,\by})  \Big)    . \nonumber
\end{align}
From Lemma 3,  we obtain
\begin{equation}  \nonumber 
\vartheta  =\delta_{P_{\br_0}} \Big(  \sum_{1 \leq j \leq 4}  m_j P_{\br_0 - \br_j+\balpha_j}  \Big)
 ,
\end{equation}
where $\br_{0} =(r_{0,1},..., r_{0,s})$, $r_{0,i} = \max_{ 1 \leq j \leq 4} r_{j,i} - \alpha_{j,i}$, $1 \leq i \leq s$. \\
Taking into account that $\max_j |m_j| \leq n^{10s}$, $\max_{i,j} \alpha_{j,i} \leq  10s \log_2 n$   and \\
 $\min_j \max_i (r_{j,i}) \geq W_0 = 50 s^2 p_0 \log_2 n$, we have
$| \sum_{1 \leq j \leq 4}  m_j P_{\br_0 - \br_j+\balpha_j} | < P_{\br_0} /10$. Hence
\begin{equation}  \nonumber 
\vartheta =
              \Delta \Big(  \sum_{1 \leq j \leq 4}  m_j P_{\br_0 - \br_j+\balpha_j} =0\Big) =\Delta \Big(  \sum_{1 \leq j \leq 4}  m_j P_{ - \br_j+\balpha_j} =0\Big).
\end{equation}
By Lemma 5, $\#U_3 \ll n^{s-1} \log^{40s} n $.  Applying Corollary 1, we obtain
\begin{equation}  \label{Le5-4}
 \sum_{ \br_j \in U_3, \; 1 \leq j \leq 4}  \Delta \Big(  \sum_{1 \leq j \leq 4}  m_j P_{ - \br_j+\balpha_j} =0\Big)
  \ll (\#U_3)^2 \ll n^{2s-2} \log_2^{80s} n .
\end{equation}
From \eqref{Le5-1} - \eqref{Le5-4} and Lemma 1,  we get
\begin{equation*} 
 \bE_{x}  \bE_{\btheta}^2 \ddot{\cD}_3^2([\theta N] )  \ll
      \sum_{ \substack{\alpha_{j,i} \in [0, 10s \log_2 n],\;
       0 < |m_j| \leq n^{10s}\\ (m_j,p_0)=1,\; j \in [1,4],\;  1 \leq i \leq s }}
       \frac{ n^{2(s-1)}  \log_2^{80s} n}{\bar{m}_1 \cdots \bar{m}_4 P_{\balpha_1 +...+\balpha_4}}
        \ll n^{2s-7/4}.
\end{equation*}
By  Chebyshev's inequality, we have
\begin{multline}  \label{Le5-6}
{\rm Prob}(\bx \in [0,1)^s \;|\; \bE_{\btheta} \ddot{\cD}_3^2([\theta N])) > n^{s-1/8})
 \leq    \bE_{x} \big(\bE_{\btheta}\ddot{\cD}_3^2([\theta N])\big)^2/ n^{2s-1/4} \\
         n^{2s-7/4}/ n^{2s-1/4}      = n^{-3/2}.
\end{multline}
Now using the Borel-Cantelli lemma, we get the assertion of Lemma 7.  \qed \\ \\
{\bf Lemma 8.}  {\it With the notations as above}
\begin{equation} \nonumber 
  \bE_{\btheta} \rD_{1,2} \ll  n^{s-1/16}  \quad \for \;a.s. \; \bx, \; s \geq 1.
\end{equation} \\
{\bf Proof.}
According to the Borel-Cantelli lemma, it is sufficient to prove that
$  {\rm Prob} ( \bx \in [0,1)^s \; | \; ( \bE_{\btheta} \rD_{1,2} )^4 > n^{4(s-1/16)}) \ll n^{-5/4} $.
%
Using  Chebyshev's inequality, we get that it is enough to verify that
 $    \bE_{x} ( \bE_{\btheta} \rD_{1,2} )^4 \ll n^{4s-3/2}$. We will prove that
 $    \bE_{x} ( \bE_{\btheta} \rD_{1,2} )^4 \ll n^{2s} <n^{4s-3/2}$ . \\

By \eqref{d-16} and \eqref{d-18},  we have that the expression of $(\bE_{\btheta} \rD_{1,2} )^4
 $ includes the summation
 over $\br_1,...,\br_8 $,
$\balpha_1,...,\balpha_8$,  $m_1,...,m_8 $ with
 $m_j/P_{\br_j-\balpha_j} = m_{j+4}/P_{\br_{j+4}-\balpha_{j+4}}$, $1 \leq j \leq 4$.
Hence $m_j=m_{j+4}$,    $\br_j - \balpha_j=\br_{j+4} - \balpha_{j+4}$, $1 \leq j \leq 4$. \\
From \eqref{d-18},  we derive
\begin{multline}  \nonumber
   \bE_{x}  ( \bE_{\btheta} \rD_{1,2} )^2 \leq
    \sum_{ \substack{\alpha_{j,i} \in [0, 10s \log_2 n]\\ 1 \leq i \leq s, 1 \leq j \leq 8}} \;
     \sum_{\substack{ \br_j \in U_1 \\  1 \leq j \leq 8}} \;
    \sum_{\substack{  0 < |m_j| \leq n^{10s} \\ (m_j,p_0)=1,1 \leq j \leq 8}} \; \prod_{j=1}^4  \Big(
     \Big| \bE_{\btheta} \big(\varphi_{ \br_j,[\theta N],m_j P_{\balpha_j}} \\
 \times  \varphi_{ \br_{j+4},[\theta N],m_{j+4} P_{\balpha_{j+4}}} \big) \Big|
 \Delta(m_j=m_{j+4}) \Delta(\br_j - \balpha_j=\br_{j+4} - \balpha_{j+4} ) \Big)  \\
   \;
 \times   \prod_{j=1}^8|\psi_{ \br_j}(m_j P_{\balpha_j},\by) |
       \Big|\bE_{x}  e\Big(\sum_{1 \leq j \leq 8} \frac{ m_j P_{\balpha_j  }}{P_{\br_j}}
       (V_{\br_j,\bx} -V_{\br_j,\by})  \Big)\Big| .
\end{multline}
Hence
\begin{multline}  \label{Le6-2}
   \bE_{x}  ( \bE_{\btheta} \rD_{1,2} )^2  \ll
    \sum_{ \substack{\alpha_{j,i} \in [0, 10s \log_2 n]\\ 1 \leq i \leq s, 1 \leq j \leq 8}} \;
     \sum_{\substack{ \br_j \in U_1 \\  1 \leq j \leq 8}} \;
    \sum_{\substack{  0 < |m_j| \leq n^{10s} \\ (m_j,p_0)=1,1 \leq j \leq 8}} \;
     \frac{1}{m_1^2 ... m_4^2 P_{\balpha_1+...+\balpha_8 }} \\
 \prod_{j=1}^4\Delta(m_j=m_{j+4}) \Delta(\br_j - \balpha_j=\br_{j+4} - \balpha_{j+4} )   \\
   \;
 \times    \Big|\bE_{x}  e\Big(\sum_{1 \leq j \leq 8} \frac{ m_j P_{\balpha_j  }}{P_{\br_j}}
       (V_{\br_j,\bx} -V_{\br_j,\by})  \Big)\Big| .
\end{multline}

Taking into account that $\max_j |m_j| \leq n^{10s}$, $\max_{i,j} \alpha_{j,i} \leq  10s \log_2 n$ and \\
 $\min_j \max_i  \; r_{j,i} \geq W_0 = 50 s^2 p_0  \log_2  n$, we have
$| \sum_{1 \leq j \leq 8}  m_j P_{\br_0 - \br_j+\balpha_j} | < P_{\br_0}/4$.
Using Lemma 3, we have
\begin{equation*} 
 \bE_{x}   e\Big(  \sum_{j=1}^8  \frac{m_j P_{\balpha_j }}{P_{\br_j }}  (V_{\br_j,\bx} -V_{\br_j,\by})
 \Big)    =
              \Delta \Big(  \sum_{j=1}^8   m_j P_{- \br_j+\balpha_j} =0\Big)
              = \Delta \Big(  \sum_{j=1}^4  2m_j P_{- \br_j+\balpha_j} =0\Big).
\end{equation*}
By Lemma 5, we get $\#U_1 \ll n^s $. Applying Corollary 1, we obtain
\begin{equation*} 
 \sum_{ \br_j \in U_1, \;1 \leq j \leq 8}
 \Delta \Big(  \sum_{1 \leq j  \leq 4}
   2m_j P_{- \br_j+\balpha_j} =0\Big)
      \prod_{j=1}^4  \Delta(\br_j - \balpha_j=\br_{j+4} - \balpha_{j+4} )
  \ll \#U_1^2 \ll n^{2s}.
\end{equation*}
%
By  \eqref{Le6-2},  we have
\begin{equation*}  
   \bE_{x}  ( \bE_{\btheta} \rD_{1,2} )^2 \leq   n^{2s}
    \sum_{ \substack{\alpha_{j,i} \in [0, 10s \log_2 n]\\ 1 \leq i \leq s, 1 \leq j \leq 8}}
       \sum_{\substack{  0 < |m_j| \leq n^{10s} \\ (m_j,p_0)=1,1 \leq j \leq 4 }}
  \frac{1}{ m^2_1  ... m_4^2 P_{\balpha_1+...+\balpha_8 }} \ll n^{2s}.
\end{equation*}
Hence Lemma 8 is proved. \qed \\ \\

{ \bf 2.4. Lower bound of the variance of $\sD_{\by,\bx}([\theta N]) )$ }.\\
In Lemma 10, we will prove that $  \bE_{\btheta} \ddot{\cD}_1^2([\theta N])  \geq \kappa_3 n^{s}  $,
with some $\kappa_3 >0$.
%
This is the main result of the section.
Lemma 9 is auxiliary. In Lemma 11 we collect the results of all previous sections.\\ \\
{\bf Lemma 9.}  {\it Let
\begin{equation*}  
   G_{\br}:=
  \sum_{\substack{  0 < |m| \leq n^{10s} \\ (m ,p_0)=1 }} \frac{1}{2\pi^2|m|^2}
 \Bigg| \sum_{ \substack{\alpha_{i} \in [0, 10s \log_2 n]\\ 1 \leq i \leq s}}
          \psi_{  \br +\balpha}(m P_{\balpha},\by) /P_{\balpha}  \Bigg|^2.
\end{equation*}
Then}
\begin{equation}  \label{Le7-0a}
\bE_{\btheta} \rD_{1,1} = \sum_{\br \in U_1} G_{\br} + O(n^{s} \log^{-15s} n).
\end{equation} \\
{\bf Proof.}  By  \eqref{d-16},  we have that if $\chi_{1,\bm,\br}=1$,
 then  $m_1=-m_2$ and\\ $ \br_2 =\br_1 -\balpha_1  + \balpha_2 $.
Applying  \eqref{d-18} and Lemma 1, we get
\begin{multline}  \label{Le7-1}
 \bE_{\btheta} \rD_{1,1}
 = \sum_{\br_1, \br_2 \in U_{1}}
    \sum_{ \substack{\alpha_{j,i} \in [0, 10s \log_2 n]\\ 1 \leq i \leq s, j=1,2}}
    \sum_{\substack{  0 < |m_j| \leq n^{10s} \\ (m_j,p_0)=1,j=1,2 }}
      \bE_{\btheta} ( \varphi_{ \br_1,[\theta N],m_1 P_{\balpha_1}} \;
      \varphi_{ \br_2,[\theta N],m_2 P_{\balpha_2}} ) \\
\times e\Big( \frac{m_1}{P_{\br_1-\balpha_1}}  (V_{\br_1,\bx} -V_{\br_1,\by})
+  \frac{m_2 }{P_{\br_2-\balpha_2}}
 (V_{\br_2,\bx} -V_{\br_2,\by}) \Big)  \psi_{ \br_1}(m_1 P_{\balpha_1},\by)  \\
%
\times   \psi_{ \br_2}(m_2 P_{\balpha_2},\by) \; \chi_{1,\bm,\br} \\
   =  \sum_{ \substack{\alpha_{j,i} \in [0, 10s \log_2 n]\\ 1 \leq i \leq s, j=1,2}}
     \;\; \sum_{\br_1, \br_2 =\br_1 -\balpha_1 +\balpha_2  \in U_1}
  \sum_{\substack{  0 < |m_1| \leq n^{10s} \\ (m_1,p_0)=1 }}
      \bE_{\btheta} ( |\varphi_{ \br_1,[\theta N],m_1 P_{\balpha_1}}|^2)  \psi_{ \br_1}(m_1 P_{\balpha_1},\by)
      \\
  \times
    \psi_{ \br_1 -\balpha_1 +\balpha_2 }(-m_1 P_{\balpha_2},\by)
  =     \sum_{ \substack{\alpha_{j,i} \in [0, 10s \log_2 n]\\ 1 \leq i \leq s, j=1,2}}
     \sum_{\br_1, \br_2=\br_1 -\balpha_1 +\balpha_2  \in U_1}
  \sum_{\substack{  0 < |m_1| \leq n^{10s} \\ (m_1,p_0)=1 }} 4\\
       \times \frac{ \bE_{\btheta}  \sin^2(2\pi m_1 [\theta N]/ P_{\br_1 - \balpha_1 } )
        }{P_{\br_1+ \br_2}|e(m_1 /P_{\br_1 - \balpha_1})-1 |^2   }
        \;  \psi_{  \br_1}(m_1 P_{\balpha_1},\by)
    \psi_{ \br_1 -\balpha_1 +\balpha_2 }(-m_1 P_{\balpha_2},\by).
\end{multline}
Taking into account that $|m_1| \leq n^{10s}$,
 $\max_{i,j} \alpha_{j,i} \leq  10s \log_2 n$  and $\min_j \max_i r_{j,i}\\ \geq W_0 = 50 s^2 p_0 \log_2 n$,
  we have
$| m_1 |^3 P_{-\br_1 + 2\balpha_1 + \alpha_2} | \ll n^{-15s}$  and
\begin{multline}  \nonumber
e(m_1 /P_{\br_1 - \balpha_1})-1 = 2\pi i m_1 /P_{\br_1 - \balpha_1} + O( (m_1 /P_{\br_1 - \balpha_1} )^2),\\
 P_{\br_1} |e(m_1 /P_{\br_1 - \balpha_1})-1| = 2\pi  |m_1| P_{\balpha_1}  + O(m_1^2 P_{-\br_1 +2 \balpha_1}),
\;\;
P_{\br_1 + \br_2} |e(m_1 /P_{\br_1 - \balpha_1})-1|^2  \\
= 4\pi^2  |m_1|^2 P_{\balpha_1 + \balpha_2}
 + O( (|m_1|^3 P_{-\br_1 +2 \balpha_1 + \balpha_2} ))=
 4\pi^2  |m_1|^2 P_{\balpha_1 +\balpha_2} +O(n^{-15s}).
\end{multline}
Hence
\begin{multline} \nonumber
  \big| \frac{1}{P_{\br_1 + \br_2} (e(m_1 /P_{\br_1 - \balpha_1})-1 )   } \big|^2   -
   \big| \frac{1}{ 4\pi^2  |m_1|^2 P_{\balpha_1 +\balpha_2}}  \big|^2 \\
    = \frac{ P_{\br_1 + \br_2} |e(m_1 /P_{\br_1 - \balpha_1})-1|^2 - 4\pi^2  |m_1|^2 P_{\balpha_1 +\balpha_2}
    }{ P_{\br_1 + \br_2} |e(m_1 /P_{\br_1 - \balpha_1})-1       |^2    4\pi^2  |m_1|^2 P_{\balpha_1+ \balpha_2}}
     \ll n^{-15s}.
\end{multline}
Using \eqref{Le2-8} and \eqref{d-10}, we have $P_{\br_1}/N =O(n^{-30s})$ for $\br_1 \in U_1$,
\begin{multline}  \nonumber
\bE_{\btheta}  \sin^2(2\pi m_1 [\theta N]/ P_{\br_1 - \balpha_1 }  ) =
1/2 +O\big( \frac{1}{ N \langle m_1 / P_{\br_1 - \balpha_1} \rangle }   \big) \\
=1/2 +O( P_{\br_1 - \balpha_1 }/N)
=1/2 +O( n^{-30s}) \quad \ad  \\
\frac{ \bE_{\btheta}  \sin^2(2\pi m_1 [\theta N]/ P_{\br_1 - \balpha_1 }  )}{
   P_{\br_1 + \br_2} |e(m_1 /P_{\br_1 - \balpha_1})-1       |^2      } =
\frac{1}{ 8\pi^2  |m_1|^2 P_{\balpha_1+ \balpha_2}}  + O(n^{-15s}) .
\end{multline}
From \eqref{Le7-1},we get
\begin{multline}  \label{Le7-5}
 \bE_{\btheta} \rD_{1,1} =
       \sum_{ \substack{\alpha_{j,i} \in [0, 10s \log_2 n]\\ 1 \leq i \leq s, j=1,2}}
     \sum_{\br_1, \br_1 -\balpha_1 +\balpha_2  \in U_1}
  \sum_{\substack{  0 < |m| \leq n^{10s} \\ (m ,p_0)=1 }} \frac{1}{2\pi^2 |m|^2 P_{\balpha_1 +\balpha_1}} \\
          \psi_{  \br_1}(m_1 P_{\balpha_1},\by)
   \; \psi_{ \br_1 -\balpha_1 +\balpha_2 }(-m_1 P_{\balpha_2},\by) + O(n^{s} \log^{-15s} n) .
\end{multline}
Bearing in mind that $|  \psi_{  \br_j+\balpha_j}(m_j P_{\balpha_j},\by)| \leq p_0 \;(j=1,2)$, we obtain
\begin{multline} \label{Le7-20}
   \sum_{ \substack{\alpha_{j,i} \in [0, 10s \log_2 n]\\ 1 \leq i \leq s, j=1,2}} \;
     \sum_{\br  \in \tilde{U}} \;
  \sum_{\substack{  0 < |m| \leq n^{10s} \\ (m ,p_0)=1 }} \; \frac{1}{2 \pi^2  |m|^2 P_{\balpha_1 +\balpha_1}} \\
      \times    \psi_{  \br+\balpha_1}(m_1 P_{\balpha_1},\by)
    \; \psi_{ \br+\balpha_2 }(-m_1 P_{\balpha_2},\by) \ll \#  \tilde{U}.
\end{multline}
Let $ \br =\br_1 -\balpha_1$. By  \eqref{d-30} and Lemma 5, we see that
\begin{multline*}
 \{\br_1 \in U_1 \; | \; \br_1 -\balpha_1 +\balpha_2  \in U_1 \} =
 \{\br \in \ZZ^s \; | \; \br+\balpha_1 \in U_1 , \; \br +\balpha_2  \in U_1 \} =U_{1_{\balpha_1,\balpha_2}}\\
= (U_1 \setminus \dot{U}_{1_{\balpha_1,\balpha_2}}) \cup
 \ddot{U}_{1_{\balpha_1,\balpha_2}}, \quad \with \quad \#\dot{U}_{1_{\balpha_1,\balpha_2}}
  +\#\ddot{U}_{1_{\balpha_1,\balpha_2}}
 \ll n^{s} \log_2^{-19s} n.
\end{multline*}
By \eqref{Le7-5} and  \eqref{Le7-20}, we have
\begin{multline*}
\bE_{\btheta} \rD_{1,1}
    =   \sum_{ \substack{\alpha_{j,i} \in [0, 10s \log_2 n]\\ 1 \leq i \leq s, j=1,2}} \;
     \sum_{\br  \in U_{1_{\balpha_1,\balpha_2}}} \;
  \sum_{\substack{  0 < |m| \leq n^{10s} \\ (m ,p_0)=1 }} \frac{1}{2 \pi^2  |m|^2 P_{\balpha_1 +\balpha_1}} \;
    \psi_{  \br+\balpha_1}(m_1 P_{\balpha_1},\by)  \\
   \times \psi_{ \br+\balpha_2 }(-m_1 P_{\balpha_2},\by) + O(n^{s} \log^{-15s} n)
  =   \sum_{ \substack{\alpha_{j,i} \in [0, 10s \log_2 n]\\ 1 \leq i \leq s, j=1,2}}
     \sum_{\br  \in U_1}
  \sum_{\substack{  0 < |m| \leq n^{10s} \\ (m ,p_0)=1 }} \frac{(2 \pi^2 )^{-1}}{|m|^2 P_{\balpha_1 +\balpha_1}} \\
      \times    \psi_{  \br+\balpha_1}(m_1 P_{\balpha_1},\by)
    \; \psi_{ \br+\balpha_2 }(-m_1 P_{\balpha_2},\by) + O(n^{s} \log^{-15s} n).
\end{multline*}
Therefore
\begin{multline}  \label{Le8-20}
\bE_{\btheta} \rD_{1,1} = \sum_{\br  \in U_1} G_{\br} + O(n^{s} \log^{-15s} n), \quad \with \quad \\
G_{\br}:=
  \sum_{\substack{  0 < |m| \leq n^{10s} \\ (m ,p_0)=1 }} \frac{1}{2 \pi^2 |m|^2}
 \Bigg| \sum_{ \substack{\alpha_{i} \in [0, 10s \log_2 n]\\ 1 \leq i \leq s}}
          \psi_{  \br +\balpha}(m P_{\balpha},\by) /P_{\balpha}  \Bigg|^2.
\end{multline}
Hence Lemma 9 is proved. \qed \\ \\
{\bf Lemma 10.}  {\it There exists $n_1 >0$ such that}
\begin{equation}  \label{Le7-0}
 \bE_{\btheta} \rD_{1,1} \geq  \kappa_3 n^s
 \quad \with \quad \kappa_3= \pi^{-2} p_0^{-6-s} 2^{-s}  \kappa_1^{2s} \kappa_2^s ,   \quad \for \quad n \geq n_1.
\end{equation} \\
{\bf Proof.}
By  Lemma 1, we get
\begin{multline}  \nonumber
 \psi_{ \br}(m,\by)  = \prod_{i=1}^s \dot{\psi}(i,\{-m M_{i,\br}/p_i \}p_i, y_{i, r_i}), \qquad \qquad
\dot{\psi}(i,m, 0) =0,  \\
 \dot{\psi}(i,m, y_{i, r_i}) = \sum_{0 \leq b  < y_{i, r_i}} e(m(b-y_{i, r_i})/p_i) \quad \for \quad y_i >0.
\end{multline}
From  \eqref{Le8-20}, we obtain
\begin{align}
&  2 \pi^2 p_0^2 (2p_0)^{2s} \geq 2 \pi^2 p_0^2G_{\br} \geq
  \sum_{\substack{  1 \leq m < p_0 \\ (m ,p_0)=1 }}
 \Bigg| \sum_{ \substack{\alpha_{i} \in [0, 10s \log_2 n]\\ 1 \leq i \leq s}}
          \psi_{  \br +\balpha}(m P_{\balpha},\by) /P_{\balpha}  \Bigg|^2 \label{Le10-80} \\
&= \sum_{ \substack{\alpha_{j,i} \in [0, 10s \log_2 n]\\ 1 \leq i \leq s, j=1,2}}
\frac{1}{P_{\balpha_1 +\balpha_2}}
 \sum_{\substack{  1 \leq m < p_0 \\ (m ,p_0)=1 }} \prod_{1 \leq i \leq s}
   \prod_{1 \leq j \leq 2}
  \dot{\psi}(i,\{(-1)^j m P_{\balpha_j} M_{i,\br+\balpha_j} /p_i \}p_i, y_{i, r_i + \alpha_{j,i}}).\nonumber
\end{align}
From \eqref{Beg-3}  and the Chinese Remainder Theorem, we get
\begin{equation*} \; \label{7-9}
       \{0 \leq m <p_0  \} = \{ \sum_{1 \leq i \leq s} m_i M_{i,\b1} p_0/p_i \mod p_0  \; | \; 0 \leq m_i < p_i, \; i \in [1,s] \},
\end{equation*}
with $\b1= (1,1,...,1)$ and $p_0 =p_1p_2 \cdots p_s$. Hence
\begin{equation*} 
       \{0 \leq m <p_0 \; | \; (m,p_0)=1  \} = \{ \sum_{1 \leq i \leq s} m_i M_{i,\b1} p_0/p_i \mod p_0 \;
        | \; 1 \leq m_i < p_i, \; i \in [1,s] \}.
\end{equation*}
Therefore
\begin{multline}  \label{Le10-40}
2 \pi^2 p_0^2 G_{\br} \geq \hslash_r :=
\sum_{ \substack{\alpha_{j,i} \in [0, 10s \log_2 n]\\ 1 \leq i \leq s, j=1,2}} \;
 \frac{1}{P_{\balpha_1 +\balpha_2}}
   \prod_{1 \leq i \leq s} \powerset_{i, \br, \balpha_1, \balpha_2},  \\
\where   \quad   \powerset_{i, \br, \balpha_1, \balpha_2}=\sum_{1 \leq m_i  < p_i}  \prod_{1 \leq j \leq 2}
    \dot{\psi}(i,\{(-1)^j m_i P_{\balpha_j} M_{i,\br+\balpha_j}  M_{i,\b1}  p_0/p_i^2 \}p_i, y_{i, r_i + \alpha_{j,i}}) \\
=\sum_{1 \leq m_i  < p_i}  \prod_{1 \leq j \leq 2} \sum_{b_j=0}^{ y_{i, r_i + \alpha_{j,i}}-1}
    e((b_j - y_{i, r_i + \alpha_{j,i}})(-1)^j m_i P_{\balpha_j} M_{i,\br+\balpha_j}  M_{i,\b1}  p_0/p_i ).
\end{multline}
We put
\begin{equation}  \label{Le9-20}
 \hslash_r = \prod_{1 \leq i \leq s}  \hslash_{r,i} \quad \with \quad
 \hslash_{r,i} =
      \sum_{ \substack{\alpha_{j,i} \in [0, 10s \log_2 n]\\  j=1,2}} \;
 \powerset_{i, \br, \balpha_1, \balpha_2}/p_i^{\balpha_{1,i} +\balpha_{2,i}}.
\end{equation}
Let $c_j = y_{i, r_i + \alpha_{j,i}}  -b_j$. Applying \eqref{Beg-1}, we have
\begin{equation} \; \nonumber
 \powerset_{i, \br, \balpha_1, \balpha_2}=
   \sum_{c_1=1}^{ y_{i, r_i + \alpha_{1,i}}}  \sum_{c_2=1}^{ y_{i, r_i + \alpha_{2,i}}}
   \Big( p_i \delta_{p_i} \big( \sum_{j=1,2} c_j(-1)^{j+1}  P_{\balpha_j} M_{i,\br+\balpha_j}
    M_{i,\b1}  p_0/p_i^2    \big)  -1   \Big).
\end{equation}
According to \eqref{Beg-3}, we have $   M_{i,\br} \equiv
	\big( P_{\br}/p_i^{r_i} \big)^{-1} \mod p_i^{r_i} $. Hence
\begin{equation} \; \nonumber
    M_{i,\br} \equiv
	\big( P_{\br}/p_i^{r_i} \big)^{-1} \mod p_i
              \quad \ad \quad
P_{\balpha_1} p_i^{-\alpha_{1,i}} M_{i,\br+\balpha_1} \equiv P_{\balpha_j} p_i^{-\alpha_{2,i}}
M_{i,\br+\balpha_2}          \mod p_i .
\end{equation}
Therefore
\begin{equation} \nonumber
 \powerset_{i, \br, \balpha_1, \balpha_2}=
   \sum_{c_1=1}^{ y_{i, r_i + \alpha_{1,i}}}  \sum_{c_2=1}^{ y_{i, r_i + \alpha_{2,i}}}
   \Big(p_i \delta_{p_i} \big( \sum_{j=1,2} c_j(-1)^{j+1}p_i^{\alpha_{j,i}}  \big)  -1   \Big).
\end{equation}
Let
\begin{multline}  \nonumber
   \hslash_{r,i}^{j_1,j_2}   =    \sum_{ \substack{\alpha_{j,i} \in [0, 10s \log_2 n]\\  j=1,2}} \;
 \powerset_{i, \br, \balpha_1, \balpha_2}^{j_1,j_2} / p_i^{\balpha_{1,i} +\balpha_{2,i}}\quad \with \quad
 \powerset_{i, \br, \balpha_1, \balpha_2}^{1,1} = \powerset_{i, \br, \balpha_1, \balpha_2}  \\
 \times   \Delta(\alpha_{1,i}>0) \Delta(\alpha_{2,i}>0),\quad \powerset_{i, \br, \balpha_1, \balpha_2}^{1,2} = \powerset_{i, \br, \balpha_1, \balpha_2}
  \Delta(\alpha_{1,i}>0) \Delta(\alpha_{2,i}=0),\\
   \powerset_{i, \br, \balpha_1, \balpha_2}^{2,1} = \powerset_{i, \br, \balpha_1, \balpha_2}
  \Delta(\alpha_{1,i}=0) \Delta(\alpha_{2,i}>0), \\
   \powerset_{i, \br, \balpha_1, \balpha_2}^{2,2} = \powerset_{i, \br, \balpha_1, \balpha_2}
  \Delta(\alpha_{1,i}=0) \Delta(\alpha_{2,i}=0).
\end{multline}
It is easy to verify that
\begin{multline}  \nonumber
\powerset_{i, \br, \balpha_1, \balpha_2}^{1,1}= (p_i-1)y_{i, r_i + \alpha_{1,i}}y_{i, r_i + \alpha_{2,i}},
\qquad \powerset_{i, \br, \balpha_1, \balpha_2}^{1,2} = -y_{i, r_i + \alpha_{1,i}}y_{i, r_i }, \\
\powerset_{i, \br, \balpha_1, \balpha_2}^{2,1} = -y_{i, r_i } y_{i, r_i + \alpha_{2,i}}
\quad \ad \quad  \powerset_{i, \br, \balpha_1, \balpha_2}^{2,2} =
 \sum_{c_1,c_2=1}^{ y_{i, r_i }}
   \big(p_i \delta_{p_i} ( c_1-c_2) -1  \big)= p_i y_{i, r_i } -y_{i, r_i }^2.
\end{multline}
By  \eqref{Le9-20}, we get
\begin{equation} \nonumber 
\hslash_{r,i} =  \sum_{j_1,j_2=1,2}  \hslash_{r,i}^{j_1,j_2},   \quad  \qquad 
 \powerset_{i, \br, \balpha_1, \balpha_2}= \sum_{j_1,j_2=1,2} 
 \powerset_{i, \br, \balpha_1, \balpha_2}^{j_1,j_2}.
\end{equation}
Hence
\begin{equation} \nonumber
 \hslash_{r,i} =  y_{i, r_i }(p_i - y_{i, r_i }) -2  y_{i, r_i } \beta_{i, r_i } +(p_i-1) \beta^2_{i, r_i }
  \; \with \;  \beta_{i, r_i }  =
  \sum_{1 \leq \alpha \leq [0, 10s \log_2 n]}  y_{i, r_i + \alpha}/p_i^{\alpha} .
\end{equation}
We consider the case $1 \leq y_{i, r_i }$ and $ \{ y_{i } p_i^{r_i} \}   \leq 1 - \kappa_1$.
If $p_i=2$, then
\begin{equation*}
 \hslash_{r,i} =  1 -2  \beta_{i, r_i } + \beta^2_{i }  =(1 - \beta_{i, r_i })^2  \geq
                    (1 - \{ y_{i, r_i } p_i^{r_i} \})^2  \geq \kappa_1^2.
\end{equation*}
If $p_i \geq 3$ and $   y_{i, r_i }  \leq p_i-2$, then
\begin{equation*}
 \hslash_{r,i} \geq  2 y_{i, r_i }  -2 y_{i, r_i } \beta_{i, r_i } + 2\beta^2_{i, r_i }
 \geq 1 - \beta_{i, r_i }  \geq
                  1 - \{ y_{i} p_i^{r_i} \}     \geq \kappa_1.
\end{equation*}
If $p_i \geq 3$ and $   y_{i, r_i }  = p_i-1$, then
\begin{equation*}
 \hslash_{r,i} =  p_i-1   -2 (p_i-1) \beta_{i, r_i } +  (p_i-1) \beta^2_{i, r_i }
  \geq (p_i-1)( 1 - \beta_{i, r_i }  )^2 \geq
                 ( 1 - \{ y_{i } p_i^{r_i} \})^2     \geq \kappa_1^2.
\end{equation*}
Therefore $ \hslash_{r,i} \geq \kappa_1^2 $ for all $p_i \geq 2$. Using \eqref{Le10-40}, we derive
\begin{equation*}
2 \pi^2 p_0^2 G_{\br} \geq \hslash_r   \geq \kappa_1^{2s} \quad \for \quad
 1 \leq y_{i, r_i } \quad \ad \quad \{ y_{i } p_i^{r_i} \}   \leq 1 - \kappa_1, \quad 1 \leq i \leq s.
\end{equation*}
Applying \eqref{In10}, we have
\begin{equation*}
     2 \pi^2 p_0^2  \sum_{\br  \in [1, n/(2p_0 )]^s} G_{\br} \geq
        \kappa_1^{2s} \kappa_2^s  (n/(2 p_0 )-1)^s \geq  p_0^{-1-s} 2^{-s} \kappa_1^{2s} \kappa_2^s n^s.
\end{equation*}
By Lemma 5, we obtain
\begin{equation*}
       U_1 \supset  [1, n/(p_0 s)]^s \setminus U_6, \quad \with \quad U_6 =[0,2W_2]^s \cup U_2,\;
       \# U_6 \ll n^s \log_2^{-10s} n.
\end{equation*}
From  \eqref{Le10-80}, we have $ \sum_{\br  \in U_6} G_{\br} \ll n^s \log_2^{-10s} n $. Thus
\begin{equation*}
      \sum_{\br  \in U_1} G_{\br}  \geq  \pi^{-2} p_0^{-5-s} 2^{-s}  \kappa_1^{2s} \kappa_2^s n^s.
\end{equation*}
Now by Lemma 9, we get the assertion of Lemma 10. \qed \\ \\
{\bf  Lemma 11.} {\it With the notations as above}
\begin{align}
& \bE_{\btheta} \big(  (\sD_{\by,\bx}([\theta N]) ) -
 \ddot{\cD}_1([\theta N]))^2 \big) \ll n^{s} \log^{-10s} n,   \label{Cor3-1} \\
& \bE_{\btheta} \sD_{\by,\bx}([\theta N]) \ll  n^s \log^{-5 s} n,   \label{Cor3-1a} \\
& \bE_{\btheta} \big(  \sD^2_{\by,\bx}([\theta N]) -
 \ddot{\cD}_1^2([\theta N]) \big) \ll n^{s} \log^{-5s} n   \label{Cor3-2} \\
&\bE_{\btheta}   \sD^2_{\by,\bx}([\theta N])
              \leq    40 p_0^{s+2} n^s,  \label{Cor3-3} \\
&\bE_{\btheta}   \sD^2_{\by,\bx}([\theta N]) \geq   0.5 \kappa_3 n^s ,
                   \label{Cor3-4}
\end{align}
 $n \geq n_2(\bx)$, where $\kappa_3= \pi^{-2} p_0^{-6-s} 2^{-s}  \kappa_1^{2s} \kappa_2^s $, with some $n_2(\bx) >0$  for a.s. $\bx$.\\ \\
{\bf Proof.}
From \eqref{d-12}, Lemma 4, Lemma 6 - Lemma 10, we obtain
\begin{align}
&\bE_{\btheta} ( \sD_{\by,\bx}([\theta N])  ) -
       \ddot{\cD}([\theta N]))^2  \ll \log_2^{2s} n,  \nonumber \\
&  0.5 \bE_{\btheta} \ddot{\cD}^2([\theta N]) -  \log_2^{2s} n  \leq    \bE_{\btheta}  \sD^2_{\by,\bx}([\theta N])
         \leq 2 \bE_{\btheta} \ddot{\cD}^2([\theta N]) + \log_2^{2s} n, \nonumber \\
&  \bE_{\btheta} \ddot{\cD}_1^2([\theta N])  \leq 32 p_0^{s+2} n^{s} , \quad
\ddot{\cD}([\theta N]) = \sum_{1 \leq j \leq 5} \ddot{\cD}_j([\theta N]), \;\;
 \ddot{\cD}^2_1([\theta N]) = \sum_{1 \leq j \leq 3} \rD_{1,j}, \;\;\;
  \nonumber \\
&    \bE_{\btheta} \ddot{\cD}_2^2([\theta N])  \ll n^{s} \log^{-10s} n , \quad
     \bE_{\btheta} \ddot{\cD}_3^2([\theta N])  \ll n^{s-1/8},
     \nonumber \\
&  |\ddot{\cD}_4([\theta N])|+
 |\ddot{\cD}_5([\theta N])| \ll \log^{21s^2} n, \qquad
 \bE_{\btheta} \ddot{\cD}_1([\theta N])  \ll 1/n, \nonumber \\
& \bE_{\btheta} \rD_{1,1} \geq  \kappa_3  n^s  , \qquad \bE_{\btheta} \rD_{1,2} \ll  n^{s-1/16}, \qquad
 \bE_{\btheta} \rD_{1,3} \ll 1,   \label{Cor3-5}
\end{align}
 $n \geq n_3(\bx)$, with some $n_3(\bx) >0$  for a.s. $\bx$.\\
\textsf{ Consider}  \eqref{Cor3-1}. From \eqref{Cor3-5}, we get
\begin{multline*}
 \bE_{\btheta} \big(  (\sD_{\by,\bx}([\theta N]) ) -
 \ddot{\cD}_1([\theta N])\big)^2    \leq  2\bE_{\btheta} \big(  (\sD_{\by,\bx}([\theta N]) ) -
 \ddot{\cD}([\theta N]) \big)^2   \\
 +    2 \bE_{\btheta} \big(   \ddot{\cD}([\theta N]) -
 \ddot{\cD}_1([\theta N])\big)^2  \ll  \log^{2 s} n +  \bE_{\btheta} \big(  \sum_{j=2}^5
 \ddot{\cD}_j([\theta N])\big)^2 \\
   \ll  \log^{2 s} n + \sum_{j=2}^5 \bE_{\btheta}  \ddot{\cD}^2_j([\theta N])
 \ll n^s \log^{-10s} n.
\end{multline*}
\textsf{ Consider}  \eqref{Cor3-1a}. By \eqref{Cor3-1} and \eqref{Cor3-5}, we have
\begin{multline*}
 |\bE_{\btheta} \big(  (\sD_{\by,\bx}([\theta N]) ) -
 \ddot{\cD}_1([\theta N]))|    \leq \big(\bE_{\btheta} \big(  (\sD_{\by,\bx}([\theta N]) ) -
 \ddot{\cD}_1([\theta N]))^2 \big)^{1/2} \\
  \ll n^s  \log^{-5 s} n, \qquad
  \bE_{\btheta} \ddot{\cD}_1([\theta N]) \ll 1/n,
\end{multline*}
and \eqref{Cor3-1a} follows. \\
\textsf{ Consider}  \eqref{Cor3-2}. From \eqref{Cor3-1} and \eqref{Cor3-5}, we derive
\begin{multline*}
\big| \bE_{\btheta} \big(  \sD^2_{\by,\bx}([\theta N]) -
 \ddot{\cD}_1^2([\theta N]) \big) \big| \leq   \bE_{\btheta}\big|\big(  \sD_{\by,\bx}([\theta N]) -
 \ddot{\cD}_1([\theta N]) \big)   \\
  \times  \big(  \sD_{\by,\bx}([\theta N]) -
 \ddot{\cD}_1([\theta N]) \big) \big| \leq
    \big(\bE_{\btheta}   (\sD_{\by,\bx}([\theta N]) -
 \ddot{\cD}_1([\theta N]))^2 \big)^{1/2}  \\
  \times    \big(\bE_{\btheta}   (\sD_{\by,\bx}([\theta N]) +
 \ddot{\cD}_1([\theta N]))^2 \big)^{1/2} \\
  \ll  n^{s/2} \log^{-5s} n
   \Big(\bE_{\btheta}   \big( (\sD_{\by,\bx}([\theta N])-\ddot{\cD}_1([\theta N])) +
 2\ddot{\cD}_1([\theta N])\big)^2 \Big)^{1/2} \\
\ll   n^{s/2} \log^{-5s} n    \Big( 2\bE_{\btheta}
  \big( \sD_{\by,\bx}([\theta N])-\ddot{\cD}_1([\theta N])\big)^2
 + 8\bE_{\btheta}   \big( \ddot{\cD}_1([\theta N])\big)^2
  \Big)^{1/2} \\
  \ll n^{s/2} \log^{-5s} n ( n^s \log^{-10s} +  n^s)^{1/2} \ll n^{s} \log^{-5s} n .
\end{multline*}
\textsf{ Consider}  \eqref{Cor3-3}. By \eqref{Cor3-2} and \eqref{Cor3-5}, we have
\begin{equation*} 
\bE_{\btheta}   \sD^2_{\by,\bx}([\theta N]) =
 \big( \bE_{\btheta} \big(  \sD^2_{\by,\bx}([\theta N]) -
 \ddot{\cD}_1^2([\theta N]) \big) \big)
  + \bE_{\btheta}   \ddot{\cD}_1^2([\theta N])  \leq 40 p_0^{s+2} n^{s},
\end{equation*}
 $n \geq n_4(\bx)$, with some $n_4(\bx) >0$  for a.s. $\bx$.\\
\textsf{ Consider}  \eqref{Cor3-4}.
Bearing in mind that $\ddot{\cD}^2_1([\theta N]) = \sum_{1 \leq j \leq 3} \rD_{1,j}([\theta N]) $,
we obtain from \eqref{Cor3-5}
\begin{equation*}
 \kappa_3  n^s \leq  \bE_{\btheta} \rD_{1,1} \leq  \bE_{\btheta} \ddot{\cD}_1^2([\theta N])
   + |\bE_{\btheta} \rD_{1,2} | + |\bE_{\btheta} \rD_{1,3} | =
   \bE_{\btheta} \ddot{\cD}_1^2([\theta N]) +O(n^{s-1/6}).
\end{equation*}
Using \eqref{Cor3-2}, we get
\begin{equation*}
  \bE_{\btheta}   \sD^2_{\by,\bx}([\theta N]) \geq
  \bE_{\btheta} \ddot{\cD}_1^2([\theta N]) -n^{s}\log_2^{-5s} n  \geq
 0.5 \kappa_3 n^s ,
\end{equation*}
 $n \geq n_5(\bx)$, with some $n_5(\bx) >0$  for a.s. $\bx$.
Hence, \eqref{Cor3-4} and  Lemma 11 are proved. \qed \\ \\
{\bf 2.5. Four moments  and L\'evy conditional variance  estimates.} \\ \\
%
{\bf Lemma 12.} {\it With the notations as above}
\begin{equation}  \label{Le8-0}
\varpi:=    \sum_{1 \leq k \leq W_1} \bE_{\btheta} \DD^4_{k}([\theta N]) \ll n^{2s -1/8} \quad \for \; a.s. \; \bx.
\end{equation} \\
{\bf Proof.} Let $\dk_{j,k_1,k_2} =k_1$ for $j \in [1,4]$ and $\dk_{j,k_1,k_2} =k_2$ for $j \in [5,8]$.\\
By Lemma 1, \eqref{Le3-1} and \eqref{d-10},  we get
\begin{multline} \label{Le8-2}
\bE_{x} \varpi^2= \sum_{1 \leq k_1,k_2 \leq W_1}  \bE_{x} \Big( \bE_{\btheta}\Big( \sum_{\br \in \UU_{k_1}}
 \ddot{\dD}_{\br,L} \Big)^4 \bE_{\btheta}\Big( \sum_{\br \in \UU_{k_2}}    \ddot{\dD}_{\br,L} \Big)^4 \Big)    \\
 = \sum_{1 \leq k_1,k_2 \leq W_1}
 \sum_{\substack{ \br_j \in \UU_{\dk_{j,k_1,k_2}} \\ 1 \leq j \leq 8}}
    \sum_{ \substack{\alpha_{j,i} \in [0, 10s \log_2 n]\\ 1 \leq i \leq s, 1 \leq j \leq 8}}
    \sum_{\substack{  0 < |m_j| \leq n^{10s} \\ (m_j,p_0)=1,1 \leq j \leq 8 }}
       \bE_{\btheta} \Big(\prod_{j=1}^4      \varphi_{ \br_j,[\theta N],m_j P_{\balpha_j}}\Big)  \\
\times
       \bE_{\btheta} \Big(\prod_{j=5}^8      \varphi_{ \br_j,[\theta N],m_j P_{\balpha_j}}\Big)
    \prod_{j=1}^8  \psi_{ \br_j}(m_j P_{\balpha_j},\by)
    \bE_{x}   e\Big(  \sum_{j=1}^8  \frac{m_j}{P_{\br_j-\balpha_j}} (V_{\br_j,\bx} -V_{\br_j,\by})
 \Big)     \\
\ll
  \sum_{ \substack{\alpha_{j,i} \in [0, 10s \log_2 n]\\ 1 \leq i \leq s, 1 \leq j \leq 8}}
    \sum_{\substack{  0 < |m_j| \leq n^{10s} \\ (m_j,p_0)=1,1 \leq j \leq 8 }}
    \prod_{1 \leq j \leq 8} \frac{1}{|m_j P_{\alpha_j}|}  \sum_{1 \leq k_1,k_2 \leq W_1}  \\
   \times
 \sum_{\substack{ \br_j \in \UU_{\dk_{j,k_1,k_2}} \\ 1 \leq j \leq 8}}
              \bE_{x}   e\Big(  \sum_{1 \leq j \leq 8}  \frac{m_j P_{\balpha_j}}{P_{\br_j}}
               (V_{\br_j,\bx} -V_{\br_j,\by})
 \Big)     .
\end{multline}
From Lemma 3,  we obtain
\begin{equation*} 
   \bE_{x}   e\Big(  \sum_{1 \leq j \leq 8}  \frac{m_j P_{\balpha_j} }{P_{\br_j}}  (V_{\br_j,\bx} -V_{\br_j,\by})
 \Big)    =
\delta_{P_{\br_0}} \Big(  \sum_{1 \leq j \leq 8}  m_j P_{\br_0 - \br_j+\balpha_j}  \Big) ,
\end{equation*}
where $ \br_0 =(r_{0,1},...,r_{0,s})$, $r_{0,i} = \max_{ 1 \leq j \leq 8} (0,r_{j,i} - \alpha_{j,i})$, $1 \leq i \leq s$. \\
Taking into account that $\max_j |m_j| \leq n^{10s}$, $\max_{i,j} \alpha_{j,i} \leq  10s \log_2 n$ \\ and
$\min_j \max_i(r_{j,i}  ) \geq  W_0 = 50 s^2 p_0 \log_2 n$, we have
$ \max_{1 \leq j \leq 8}|   m_j P_{\br_0 - \br_j+\balpha_j} | < P_{\br_0} /16$.
By Lemma 3
\begin{equation*}  
 \bE_{x}   e\Big(  \sum_{1 \leq j \leq 8}  \frac{m_j P_{\balpha_j}}{P_{\br_j}}  (V_{\br_j,\bx} -V_{\br_j,\by})
 \Big)    =
              \Delta \Big(  \sum_{1 \leq j \leq 8}  m_j P_{- \br_j+\balpha_j} =0\Big).
\end{equation*}
Applying Corollary 1 and Lemma 5, we have
\begin{equation}  \label{Le8-5}
 \sum_{\substack{ \br_j \in \UU_{\dk_{j,k_1,k_2}} \\ 1 \leq j \leq 8}}
   \Delta \Big(  \sum_{1 \leq j \leq 8}  m_j P_{ - \br_j+\balpha_j} =0\Big)
  \ll \max_{k}  (\# \UU_k)^4 \ll n^{4(s-1)} \log_2^{80s} n.
\end{equation}
From \eqref{Le8-2}-\eqref{Le8-5} and \eqref{d-10},  we get
\begin{multline*}  
    \bE_{x} \varpi^2 \ll    \sum_{ \substack{\alpha_{j,i} \in [0, 10s \log_2 n]\\ 1 \leq i \leq s, 1 \leq j \leq 8}}
    \sum_{\substack{  0 < |m_j| \leq n^{10s} \\ (m_j,p_0)=1,1 \leq j \leq 8 }}
    \prod_{1 \leq j \leq 8} \frac{1}{|m_j P_{\alpha_j}|} \sum_{1 \leq k_1,k_2 \leq W_1} n^{4(s-1)} \log_2^{80s} n \\
 \ll      \log_2^{8}n \;
       W_1^2  n^{4(s-1+0.05)} \ll \log_2^{8}n
       (n/ \log^{20s}_2 n)^2  n^{4s-4+0.2)}   \ll n^{4s-3/2}.
\end{multline*}
By  Chebyshev's inequality, we have
\begin{equation*} 
P(  \varpi > n^{2s-1/8})  \leq    \bE_{x} \varpi^2/ n^{4s-1/4} \ll n^{4s-3/2 -4s +1/4}=n^{-5/4}.
\end{equation*}
Now using the Borel-Cantelli lemma, we get the assertion of Lemma 12. \qed  \\ \\
Denote by $\dot{\cF}(l)$ the sigma field on $[0,1)$ generated by
  $\{ [\frac{j}{2^l},\frac{j+1}{2^l}) \; | \;
   j=0,..., 2^l -1  \}$. Let $l(0)=0, l_k = (k+1) W_2+W_3$, $\cF_k=\dot{\cF}(l_k)$.
 We will consider the probability space $([0,1), B([0,1)), \lambda_1)$ and the  conditional expectation
\begin{equation}\label{Le8-12}
  {\bf E}_{\btheta} ( f(\theta) \; | \; \cF_{k})  =
       2^{l_k} \int_{\fn/ 2^{l_{k}}}^{(\fn+1)/ 2^{l_{k}}}  f(x )    dx \qquad \for \quad
    \theta \in   [\frac{\fn}{2^{l_{k}}}, \frac{\fn +1}{2^{l_{k}}}].
\end{equation}  \\
{\bf Lemma 13.} {\it With the notations as above}
\begin{multline} \label{Le12-10}
  {\bf E} (  \DD_{k}([\theta N]) \; | \; \cF_{k-1})  =O(n^{-30s}),
  \;\;\;       {\bf E} (  \DD_{k}([\theta N]) \; | \; \cF_{k})  =  \DD_{k}([\theta N]) + O(n^{-30s}) ,  \\
   {\bf E}_{\btheta} (  \DD^2_{k}([\theta N]) \; | \; \cF_{k-1} ) =
          {\bf E}_{\btheta} (  \DD^2_{k}([\theta N]) ) +O(n^{-30s}).
\end{multline}  \\
{\bf Proof.}
By  Lemma 1, \eqref{Le3-1}, \eqref{d-10} and  \eqref{Le9-4},  we get
\begin{multline*}  
   {\bf E}_{\btheta} (  \DD_{k}([\theta N]) \; | \; \cF_{k-1}) =  \sum_{ \br \in \UU_k}
   \sum_{ \substack{\alpha_i \in [0, 10s \log_2 n]\\ 1 \leq i \leq s}}
    \sum_{\substack{ 0 < |m| \leq n^{10s} \\ (m,p_0)=1 }}
   \frac{\sqrt{-1} }{P_{\br}(e(m/P_{\br-\balpha})-1)}
      \\
  \times 2{\bf E}_{\btheta} ( \sin(2\pi m  [\theta N] /P_{\br-\balpha} )  \; | \; \cF_{k-1})
        \psi_{ \br}(m P_{\balpha},\by)\;
e\Big( \frac{m P_{\balpha}}{P_{\br}}  (V_{\br,\bx} -V_{\br,\by}) \Big) \\
  \ll n^{-40s} \sum_{ \br \in \UU_k}
   \sum_{ \substack{\alpha_i \in [0, 10s \log_2 n]\\ 1 \leq i \leq s}}
    \sum_{\substack{ 0 < |m| \leq n^{10s} \\ (m,p_0)=1 }}
   \frac{1 }{\bar{m} P_{\balpha}}     \ll n^{-30s}.
\end{multline*}
Similarly, using \eqref{Le9-3}, we have
\begin{multline*}  
   {\bf E}_{\btheta} (  \DD_{k}([\theta N]) \; | \; \cF_{k}) - \DD_{k}([\theta N])
    =  \sum_{ \br \in \UU_k}
   \sum_{ \substack{\alpha_i \in [0, 10s \log_2 n]\\ 1 \leq i \leq s}}
    \sum_{\substack{ 0 < |m| \leq n^{10s} \\ (m,p_0)=1 }}
   \frac{\sqrt{-1}}{P_{\br}(e(m/P_{\br-\balpha})-1)}
      \\
  \times    2\Big( {\bf E}_{\btheta} \big( \sin(2\pi m  [\theta N] /P_{\br-\balpha} )  \; | \; \cF_{k}\big) -
    \sin(2\pi m  [\theta N] /P_{\br-\balpha} )  \Big)    \psi_{ \br}(m P_{\balpha},\by)
 \\
  \times e\Big( \frac{m P_{\balpha}}{P_{\br}}  (V_{\br,\bx} -V_{\br,\by}) \Big) \ll n^{-40s} \sum_{ \br \in \UU_k}
   \sum_{ \substack{\alpha_i \in [0, 10s \log_2 n]\\ 1 \leq i \leq s}}
    \sum_{\substack{ 0 < |m| \leq n^{10s} \\ (m,p_0)=1 }}
   \frac{1 }{\bar{m}P_{\balpha}}    \ll n^{-30s}.
\end{multline*}
Therefore, the first part of Lemma 13 is proved.\\
\textsf{ Consider} \eqref{Le12-10}.    By \eqref{Le3-1} and \eqref{d-10},  we obtain
\begin{multline*} 
 \DD_{k}^2([\theta N]) = - \sum_{\br_1, \br_2 \in \UU_{k}}
    \sum_{ \substack{\alpha_{j,i} \in [0, 10s \log_2 n]\\ 1 \leq i \leq s,j=1,2}}
    \sum_{\substack{  0 < |m_j| \leq n^{10s} \\ (m_j,p_0)=1,j=1,2 }} 4 \\
  \times   \frac{\sin(2\pi m_1  [\theta N] /P_{\br_1 -\balpha_1} ) \sin(2\pi m_2  [\theta N] /P_{\br_2 -\balpha_2} ) }{P_{\br_1}(e(m_1/P_{\br_1-\balpha_1})-1) P_{\br_2}(e(m_2/P_{\br_2-\balpha_2})-1)}
      \\
  \times
    \psi_{ \br_1}(m_1 P_{\balpha_1},\by)  \psi_{ \br_2}(m_2 P_{\balpha_2},\by)
   \; \;
e\Big( \frac{m_1 P_{\balpha_1}}{P_{\br_1}}  (V_{\br_1,\bx} -V_{\br_1,\by})  +  \frac{m_2 P_{\balpha_2}}{P_{\br_2}}
 (V_{\br_2,\bx} -V_{\br_2,\by}) \Big) .
\end{multline*}
Now using  \eqref{Le9-5},  we get
\begin{multline*}  
 \Big| \bE_{\btheta} (  \DD_{k}^2([\theta N])  | \cF_{k-1} ) -\bE_{\btheta}   \DD_{k}^2([\theta N])        \Big|  \ll
  \sum_{\br_1, \br_2 \in \UU_{k}}
    \sum_{ \substack{\alpha_{j,i} \in [0, 10s \log_2 n]\\ 1 \leq i \leq s,j=1,2}}
    \sum_{\substack{  0 < |m_j| \leq n^{10s} \\ (m_j,p_0)=1,j=1,2 }} 1 \\
  \times   \Big| \bE_{\btheta} \Big( \sin(2\pi m_1  [\theta N] /P_{\br_1 -\balpha_1} )
  \sin(2\pi m_2  [\theta N] /P_{\br_2 -\balpha_2} )  \; |\; \cF_{k-1} \Big) \\
   -
  \bE_{\btheta} \Big( \sin(2\pi m_1  [\theta N] /P_{\br_1 -\balpha_1} )
   \sin(2\pi m_2  [\theta N] /P_{\br_2 -\balpha_2} )    \Big)\Big|
      \frac{1}{|m_1 m_2|P_{\balpha_1 +\balpha_2} } \\
  \ll n^{-40s} \sum_{\br_1, \br_2 \in \UU_{k}}
    \sum_{ \substack{\alpha_{j,i} \in [0, 10s \log_2 n]\\ 1 \leq i \leq s,j=1,2}}
    \sum_{\substack{  0 < |m_j| \leq n^{10s}
    \\ (m_j,p_0)=1,j=1,2 }} \frac{1}{|m_1 m_2|P_{\balpha_1 +\balpha_2} } \ll  n^{-30s}.
\end{multline*}
Hence Lemma 13 is proved. \qed \\ \\

{\bf 2.6. Martingale approximation.}\\
Let
\begin{equation}\label{Le12-14}
      \xi_k = {\bf E} [  \DD_{k}([\theta N]) \; | \; \cF_{k}]
      - {\bf E} [  \DD_{k}([\theta N]) \; | \; \cF_{k-1}], \; k=1,2,...
\end{equation}
Then $(\xi_k)_{k \geq 1}$ is the martingale difference array satisfying
${\bf E} [ \xi_k  |  \cF_{k-1}]=0, \; k=1,2,...$ .
Bearing in mind  \eqref{d-10} and \eqref{d-12}, we define
\begin{equation}\label{Le13-40}
 \dot{\SS}_n: = \sum_{k=1 }^{ W_1} \xi_k,  \quad   \ddot{\SS}_n:
   =\sum_{k=1 }^{ W_1} \DD_{k}([\theta N])=\ddot{\cD}_1 ([\theta N]), \;\;\; \dot{\varrho}^2_n:=   {\bf E}_{\btheta}(\dot{\SS}^2_n),
   \;\; \; \ddot{\varrho}^2_n:=   {\bf E}_{\btheta}(\ddot{\SS}^2_n).
\end{equation} \\
{\bf Lemma 14.} {\it With the notations as above}
\begin{equation} \nonumber 
           \DD_{k}([\theta N]) - \xi_k = O(n^{-30s}) , \quad
           \ddot{\SS}_n - \dot{\SS}_n = O(n^{-29s}), \quad \DD_{k}([\theta N])^2 - \xi_k^2 = O(n^{-28s}),
\end{equation}
\begin{equation*}  
   \ddot{\varrho}^2_n - \dot{\varrho}^2_n = O(n^{-11s})      \qquad \quad \ad
   \qquad  \quad  |\xi_k|^4 \leq  8| \DD_{k}([\theta N])|^4+ O(n^{-30s}) .
\end{equation*}\\
{\bf Proof.}
By \eqref{Le12-14} and Lemma 13,  we get
\begin{multline*}
 \DD_{k}([\theta N]) - {\bE}_{\theta} [  \DD_{k}([\theta N]) \; | \; \cF_{k}] \ll n^{-30s}, \qquad
\;\; {\bE}_{\theta} ( \DD_{k}([\theta N]) \; | \; \cF_{k-1}) \ll n^{-30s}  \\
 \DD_{k}([\theta N]) - \xi_k = \DD_{k}([\theta N]) - {\bf E}_{\btheta} [  \DD_{k}([\theta N]) \; | \; \cF_{k}]
 -  {\bf E}_{\btheta} ( \DD_{k}([\theta N]) \; | \; \cF_{k-1})  \ll n^{-30s},
\end{multline*}
$ \ddot{\SS}_n - \dot{\SS}_n = O(n^{-29s})$. From \eqref{d-11} and  \eqref{Le13-40},
  we have $\DD_k([\theta N]) \ll n^{s+1}$, and
 $|\dot{\SS}_n|  + |\ddot{\SS}_n| \ll n^{s+2}$.
Hence
%
\begin{equation*}
 | \DD_{k}([\theta N])^2 - \xi_k^2| \leq  |  \DD_{k}([\theta N]) - \xi_k| \;
( 2| \DD_{k}([\theta N])| +  | \DD_{k}([\theta N]) - \xi_k|)
     \ll n^{-29s+1}
\end{equation*}
 and
\begin{multline*}
 | \ddot{\varrho}^2_n - \dot{\varrho}^2_n| = |{\bE}_{\theta} (\ddot{\SS}_n^2)-
  {\bE}_{\theta}(\dot{\SS}_n^2) |    =  \big|{\bE}_{\theta} \big((\ddot{\SS}_n - \dot{\SS}_n)(\ddot{\SS}_n + \dot{\SS}_n)\big) \big|
  \\
    \leq  \big({\bE}_{\theta} \big((\ddot{\SS}_n - \dot{\SS}_n)^2 \big)
    {\bE}_{\theta} \big((\ddot{\SS}_n + \dot{\SS}_n)^2\big) \Big)^{1/2}
    \ll n^{(-29s+2s +4)/2}     \ll n^{-11s}.
\end{multline*}
By H\"older's inequality, we have $(a +b)^4 \leq 8(a^4 + b^4)$. Therefore
\begin{multline*}
   |\xi_i|^4 = |\DD_{k}([\theta N]) + \xi_k -\DD_{k}([\theta N])|^4 \leq 8| \DD_{k}([\theta N])|^4 +8  | \DD_{k}([\theta N]) - \xi_i|^4  \\
   =  8| \DD_{k}([\theta N])|^4 +O(n^{-40s}).
\end{multline*}
Hence Lemma 14 is proved. \qed
\\

 We shall use  the following variant of the {\it martingale central limit theorem} (see [Ha, p.~58,
 Corollary 2.1]):

 Let $(\Omega, \cF,P)$ be a probability space and $\{(\zeta_{n,k}, F_{n,k}) \;| \;  k=1,...,\ell_n\}$ be a martingale difference array with
 ${\bf E} [\zeta_{n,k} | F_{n,k-1}] =0$ a.s. ($F_{n,0}$ is the trivial field).
 \\ \\
{\bf Theorem  D.} {\it Let
\begin{multline}\label{Le13-10}
 \SS_{n} =\sum_{1 \leq k \leq \ell_n} \zeta_{n,k}, \quad \quad
 L(n,\epsilon) = \sum_{1 \leq k \leq \ell_n} {\bf E} (\zeta_{n,k}^2 \delta(|\zeta_{n,k}|>\epsilon)),
  \quad \quad
  \sum_{1 \leq k \leq \ell_n}  {\bf E} (\zeta_{n,k}^2) =1,\\
  \VV_{n}^2 =\sum_{1 \leq k \leq \ell_n} {\bf E} (\zeta_{n,k}^2 | F_{n,k-1}) ,  \quad
  L(n,\epsilon)  \stackrel{P}{\rightarrow} 0
    \quad \forall  \epsilon>0,  \quad   \VV_{n}^2 \stackrel{P}{\rightarrow} 1.
\end{multline}
Then $ \SS_{n} 	\stackrel{w}{\rightarrow} \cN(0,1)    $. } \\

   Now we apply Theorem D to  the martingale difference array $\{(\zeta_{n,k}, F_{n,k}) \;| \;
     k=1,...,\ell_n\}$ with $F_{n,k}= \cF_k$,
 $\zeta_{n,k} = \xi_k/\dot{\varrho}_n$
 and  $\ell_n= W_1$.\\ \\
{\bf Lemma 15.} {\it With the notations as above}
   $ \dot{\SS}_n / \dot{\varrho}_n 	\stackrel{w}{\rightarrow} \cN(0,1)  $. \\ \\
{\bf Proof.}
  By (\ref{Le12-14}),  $(\xi_i)_{i \geq 1}$ is the martingale difference sequence (and consequently orthogonal).
From (\ref{Cor3-2}), \eqref{Cor3-4}, (\ref{Le13-40}) and Lemma 14, we obtain
\begin{multline}  \label{Le14-1}
\sum_{k \in [1,W_1]}  {\bf E}_{\btheta} \xi_k^2  =
 {\bf E}_{\btheta}  \Big(\sum_{k \in [1,W_1]} \xi_k  \Big)^2 =  \dot{\varrho}_n^2
 =  \ddot{\varrho}_n^2 +O(n^{-11s})
    =    {\bf E}_{\btheta} \Big(\sum_{k \in [1,W_1]}  \DD_{k}([\theta N]) \Big)^2 \\
     +O(n^{-11s}) = {\bf E}_{\btheta}  \ddot{\cD}^2_{1}([\theta N])
              +O(n^{-11s})
                 = {\bf E}_{\btheta} (  \sD^2_{\by,\bx}([\theta N])  ) \\
                +O(n^{s} \log_2^{-5s})  \in n^{s} [0.4 \kappa_3,\omega_3],
\end{multline}
with some $\omega_3>0$.\\
 We derive from  (\ref{Le13-10}), Lemma 12 and Lemma 14 that
\begin{multline}\nonumber
  L(n,\epsilon) = \sum_{1 \leq k \leq W_1}
   \int_0^1   |\xi_k / \dot{\varrho}_n|^{2} \delta(\xi_k/\dot{\varrho}_n|> \epsilon)  d \theta
   \leq \sum_{1 \leq k \leq W_1}\epsilon^{-2}    \int_0^1  |\xi_k / \dot{\varrho}_n|^{4} d \theta \\
   \ll \sum_{1 \leq k \leq W_1}\epsilon^{-2}  n^{-2s}
   ({\bf E}_{\btheta} \DD^4_{k,[\theta N]} +O(n^{-30s}) )
   \ll \epsilon^{-2}  n^{-2s} n^{2s-1/8} \ll n^{-1/8} .
\end{multline}
By Lemma 14, we get $  \DD_{k}([\theta N]) - \xi_k = O(n^{-10s}) $.
Using (\ref{Le14-1}) and  Lemma 13, we derive
\begin{multline*} 
   \VV_n^2 -1  = \dot{\varrho}_n^{-2} \sum_{1 \leq k \leq W_1}
   \Big( {\bf E}_{\btheta} (\xi^2_k | \cF_{k-1})
    -  {\bf E}_{\btheta} \xi_k^2 \Big)  \\
         =      \dot{\varrho}_n^{-2} \sum_{1 \leq k \leq W_1}  \Big( {\bf E}_{\btheta} (\DD^2_{k,[\theta N]}  | \cF_{k-1})
   -  {\bf E}_{\btheta} \DD^2_{k,[\theta N]} \Big) +O(n^{-29s}) = O(n^{-29s}).
\end{multline*}
Applying (\ref{Le14-1}) and Theorem D, we obtain the assertion of Lemma 15. \qed \\ \\
We need the following ``Converging Together Lemma'' : \\ \\
{\bf Lemma A.} [Du, p.105, ex.3.2.13] {\it If $X_n \stackrel{w}{\rightarrow} X$ and
$Z_n - X_n \stackrel{w}{\rightarrow} 0$,
then $ Z_n \stackrel{w}{\rightarrow}X$.}\\ \\
{ \bf 2.7 The end of the proof of the Theorem.} Let
\begin{equation} \label{LemA}
\dddot{\SS}_n =\sD_{\by,\bx}([\theta N])  \qquad \ad \qquad
     \dddot{\varrho}_n =(\bE_{\btheta}\dddot{\SS}^2_n    )^{1/2} .
\end{equation}
By  \eqref{Le13-40}, Lemma 11 and  Lemma 14, we obtain
\begin{equation} \nonumber
   \dot{\SS}_n - \ddot{\SS}_n \ll n^{-29s}, \;
 \dot{\varrho}^2_n - \ddot{\varrho}^2_n \ll n^{-11s}, \;
 \bE_{\btheta}(\ddot{\SS}_n - \dddot{\SS}_n)^2 \ll n^{s} \log_2^{-10s},   \;\;
  \ddot{\varrho}^2_n - \dddot{\varrho}^2_n \ll n^{s} \log_2^{-10s}.
\end{equation}
According to  \eqref{Cor3-4} and \eqref{Le14-1}, we get $ \dot{\varrho}^2_n \geq 0.4 \kappa_3 n^s $.
Hence
\begin{multline*} \nonumber 
  \dot{\varrho}^2_n - \dddot{\varrho}^2_n =(\dot{\varrho}^2_n - \ddot{\varrho}^2_n)
+(\ddot{\varrho}^2_n - \dddot{\varrho}^2_n) \ll n^{s} \log_2^{-5s} n, \qquad
\dddot{\varrho}^2_n \geq
 \pi^{-2} p_0^{-8-s} 2^{-s}  \kappa_1^{2s} \kappa_2^s n^s \\
     \ad \;\;\;\; \bE_{\btheta} (\dot{\SS}_n -\dddot{\SS}_n)^2
    \leq 2{\bf E_{\btheta}} \big( (\dot{\SS}_n -\ddot{\SS}_n)^2 \big) +2{\bf E_{\btheta}}
    \big( (\ddot{\SS}_n -\dddot{\SS}_n)^2 \big) \ll n^{s} \log_2^{-10s} n.
\end{multline*}
 Therefore
\begin{equation} \nonumber 
    | \frac{1}{\dot{\varrho}_n} -   \frac{1}{\dddot{\varrho}_n} |
      =
        \frac{ | \dot{\varrho}_n - \dddot{\varrho}_n |}{\dot{\varrho}_n  \dddot{\varrho}_n}
     =
       \frac{ |\dot{\varrho}_n^2- \dddot{\varrho}_n^2| } {\dot{\varrho}_n
       \dddot{\varrho}_n (\dot{\varrho}_n+ \dddot{\varrho}_n)} \ll n^{s/2} \log_2^{-5s} n.
\end{equation}
%
Applying  Lemma 11, we derive
\begin{multline*}
 \bE_{\btheta}( \frac{ \dot{\SS}_n}{\dot{\varrho}_n}  - \frac{ \dddot{\SS}_n}{\dddot{\varrho}_n} )^2 =
 \bE_{\btheta}( \frac{ \dot{\SS}_n -\dddot{\SS}_n }{\dot{\varrho}_n} +
  \dddot{\SS}_n ( \frac{1}{\dot{\varrho}_n}  -\frac{1}{\dddot{\varrho}_n})^2
     \leq 2 \dot{\varrho}^{-2}_n  \bE_{\btheta}( \dot{\SS}_n -\dddot{\SS}_n )^2  \\
     +2  (1/\dot{\varrho}_n  -1/\dddot{\varrho}_n)^2 \bE_{\btheta} \dddot{\SS}^2_n
     \ll n^{-s +s} \log_2^{-10s} n +
     ( n^{-s} \log_2^{-10s} n ) n^s  \ll  \log_2^{-10s} n.
\end{multline*}
Hence
$\dot{\SS}_n /\dot{\varrho}_n - \dddot{\SS}_n /\ddot{\varrho}_n	\stackrel{w}{\rightarrow} 0$. Bearing in mind that
$ \dot{\SS}_n /\dot{\varrho}_n	\stackrel{w}{\rightarrow} \cN(0,1)  $ and
  Lemma~A, we get that
$\dddot{\SS}_n /\dddot{\varrho}_n	\stackrel{w}{\rightarrow} \cN(0,1)  $.
 By Lemma 11 and \eqref{LemA}, we have
 the assertion of the Theorem. \qed \\ \\

\qquad \qquad \qquad \qquad \qquad \qquad {\bf Appendix} \\ \\
{\bf Lemma 16.} {\it Let $  \alpha_{i} \in [0, 10s \log_2 n], \; 1 \leq i \leq s$ ,
$ 1 \leq |m| \leq n^{10s} $, $\br \in \UU_k$, $l_k = (k+1) W_2+W_3$, $k \in [1,W_1]$. Then }
\begin{align}
& | {\bf E}_{\btheta} (  \sin(2\pi \beta [\theta N]+\eta) \; | \; \cF_{k}) |
  \leq    \min \Big( 2,  \frac{2^{l_k}-1 }{N \llangle \beta \rrangle}\Big),   \label{Le9-0} \\
&|\bE_{\btheta}(\sin^2 (2\pi [\theta N] \beta)  \; | \; \cF_k  ) -1/2|
 \leq \min\big(1, \frac{2^{l_k -2} }{N\llangle 2\beta \rrangle} \big), \label{Le9-1} \\
& |\bE_{\btheta}(\sin (2\pi [\theta N] \beta) \sin (2\pi [\theta N] \gamma)  \; | \; \cF_k)|
     \leq  \sum_{i=1,2}
   \min\big(1, \frac{2^{l_k -2} }{N\llangle \beta +(-1)^i \gamma\rrangle} \big),  \label{Le9-2} \\
& \big|  {\bf E}_{\btheta} \big(\sin(2\pi  m[\theta N]/P_{\br -\balpha}+\eta)
\; | \; \cF_{k} \big) -\sin(2\pi   m[\theta N]/P_{\br -\balpha} +\eta ) \big| \leq n^{-40s},     \label{Le9-3} \\
%
& \big| {\bf E}_{\btheta} \big(\sin(2\pi  m[\theta N]/P_{\br -\balpha}
\; | \; \cF_{k-1} )\big) \big| \leq n^{-40s},     \label{Le9-4} \\
& \nonumber \\
&\Big|\bE_{\btheta}\Big(\sin \big(\frac{2\pi m_1[\theta N]}{P_{\br_1 - \balpha_1}}   \big)
 \sin \big(\frac{2\pi m_2[\theta N]}{P_{\br_2 - \balpha_2}}   \big)
  \; | \; \cF_{k-1}\Big)  \nonumber\\
&  \qquad \qquad \qquad \quad -
  \bE_{\btheta}\Big(\sin \big(\frac{2\pi m_1[\theta N]}{P_{\br_1 - \balpha_1}}   \big)
 \sin \big(\frac{2\pi m_2[\theta N]}{P_{\br_2 - \balpha_2}}   \big)\Big) \Big|
  \ll    n^{-40s}. \label{Le9-5}
%
%
\end{align} \\
{\bf Proof.} \textsf{ Consider}  \eqref{Le9-0}.
Let $\theta \in   [\frac{\fn}{2^{\ell_{k}}}, \frac{\fn +1}{2^{\ell_{k}}}]$. From \eqref{Le8-12}, we see
\begin{multline}  \nonumber 
 {\bf E}_{\btheta} (  e(\beta [\theta N]+\eta) \; | \; \cF_{k} ) =
   2^{l_k} \int_{ \fn/ 2^{l_{k}}}^{( \fn+1)/ 2^{l_{k}}}   e(\beta [x N]+\eta)    d x  \\
= \int_{0}^1   e(\beta [(N \fn +z N)/2^{l_{k}}]+\eta)    d z = \frac{1}{N} \sum_{\ell=0}^{N-1}
  e(\beta [(N \fn +\ell)/2^{l_k}]+\eta) .
\end{multline}
Let $N= N_1 2^{l_k} +N_2$, $\ell= \ell_1 2^{l_k} +\ell_2$ with $N_2,\ell_2 \in [0, 2^{l_k} )$. Hence
$ 0 \leq \ell_1 2^{l_k} \leq N-1 -\ell_2 $, $0 \leq \ell_1  \leq [(N_1 2^{l_k} +N_2-1 -\ell_2 )/2^{l_k}]=
$$ 0 \leq \ell_1 \leq  N_1 +N_3$, with $N_3 =[ (N_2 -1-\ell_2)   / 2^{l_k}  ] \in \{-1,0\}$.
Using \eqref{Beg-2a}, we derive
\begin{multline}  
 | {\bf E}_{\btheta} (  e(\beta [\theta N]+\eta) \; | \; \cF_k) |
=     \frac{1}{N} \Big| \sum_{ 0 \leq \ell_2 \leq 2^{l_k} -1}
 \sum_{0 \leq \ell_1 \leq  N_1  +N_3}
  e(\beta \ell_1 + \beta  [(N \fn +\ell_2)/2^{l_k}])    \Big| \\
 \leq   \frac{2^{l_k}}{N} \max_{ 0 \leq \ell_2 \leq 2^{l_k} -1} \Big|
 \sum_{0 \leq \ell_1 \leq  N_1  + N_3}
  e(\beta \ell_1 )   \Big|  \leq \\
 \leq   \frac{2^{l_k}}{N} \min \Big( N_1+1,  \frac{1}{2\llangle \beta \rrangle}\Big)
 \leq    \min \Big( 2,  \frac{2^{l_k} }{2N \llangle \beta \rrangle}\Big) ,
\end{multline}
and \eqref{Le9-0} follow.\\
\textsf{ Consider}  \eqref{Le9-1} and  \eqref{Le9-2}.
Bearing in mind  \eqref{Le2-8} and \eqref{Le9-0},  we obtain
\begin{multline} \nonumber 
 2|\bE_{\btheta}(\sin^2 (2\pi [\theta N] \beta)  \; | \; \cF_k  ) -1/2| =
|\bE_{\btheta}(-\cos (4\pi[\theta N] \beta)  \; | \; \cF_k  )|\nonumber \\
 = |\bE_{\btheta}(\sin (4\pi [\theta N] \beta +\pi/2)  \; | \; \cF_k  )|   \leq \min\big(2, \frac{2^{l_k} }{2N\llangle 2\beta \rrangle} \big),
\end{multline}
and
\begin{multline}
   2 |\bE_{\btheta}(\sin (2\pi [\theta N] \beta) \sin (2\pi [\theta N] \gamma)  \; | \; \cF_k)| \nonumber\\
   =  |\bE_{\btheta}(\sin (2\pi [\theta N] (\beta-\gamma +\pi/2  )- \sin (2\pi [\theta N]
    \beta +\gamma+ \pi/2)  \; | \; \cF_k)| \nonumber \\
     \leq  \sum_{i=1,2}
   \min\big(2, \frac{2^{l_k} }{2N\llangle \beta +(-1)^i \gamma\rrangle} \big). \nonumber
\end{multline}
%
\textsf{ Consider}  \eqref{Le9-3}.
By \eqref{d-10}, we have  $ W_2= [\log_2^{20s} n]$, $W_3= [\log_2^{10s} n]$,
$   A_k = n- (k+1)W_2 $, $B_k =A_k + W_2 - 2W_3$, $h(\br)  =r_1\log_2(p_1)+ \cdots + r_s \log_2(p_s)$,
 $ \UU_{k}=  \{\br \in U \; | \;  h(\br) \in [A_k,B_k) \}$. Taking into account that
 $l_k = (k+1) W_2+W_3$, we have for   $\br \in \UU_k$ :
\begin{multline} \label{Le9-10}
  \frac{P_{\br-\balpha} 2^{l_{k-1}}}{ N } \ll 2^{n -(k+1) W_2+W_2 -2W_3 -n +k W_2+W_3-1} = 2^{-W_3-1}
    \leq 2^{-\log_2^{10s} n } \leq  n^{-40s},  \\ 
 \frac{m N}{P_{\br-\balpha}} \ll 2^{(k+1)W_2 + 20s p_0 \log_2 n} \; \ad \;
 \frac{m N }{P_{\br-\balpha} 2^{l_{k}}} \ll 2^{-W_3 +20s p_0 \log_2 n} \\
 \ll 2^{-W_3/2} \ll 2^{-0.5 \log_2^{10s} n}  \ll n^{-40s}.
\end{multline}
Let   $ \theta \in [\frac{\fn}{2^{l_{k}}}, \frac{\fn +1}{2^{l_{k}}}]$, $\theta =(\fn+\theta_1)/ 2^{l_{k}}$ with $\theta_1 \in [0,1]$. We see
\begin{multline} \nonumber 
 \frac{m [\theta N]}{P_{\br-\balpha}} =   \frac{m }{P_{\br-\balpha}}
 \Big( \Big[   \frac{ \fn N}{ 2^{l_{k}} } \Big]  +   \Big[   \frac{ \theta_1 N}{ 2^{l_{k}} } \Big]
 +\epsilon    \Big) =  \frac{m }{P_{\br-\balpha}}  \Big[   \frac{ \fn N}{ 2^{l_{k}} } \Big] +O(\frac{mN }{P_{\br-\balpha} 2^{l_k}} ) \\
 =
 \frac{m }{P_{\br-\balpha}}  \Big[   \frac{ \fn N}{ 2^{l_{k}} } \Big] +O(n^{-100s}), \quad
   \with \quad\epsilon \in [-2,2].
\end{multline}
Bearing in mind that $ |\sin(x_1) -\sin(x_2)| = |x_1-x_2| |\cos(x_1 +z(x_2 -x_1))| \leq |x_1-x_2|$
 $(z \in [0,1])$, we get
\begin{equation}  
   \sin  \Big(  \frac{2\pi m [\theta N]}{P_{\br-\balpha}} +\eta \Big) -
    \sin \Big(  \frac{2\pi m }{P_{\br-\balpha}}  \Big[   \frac{ \fn N}{ 2^{l_{k}} } \Big] +\eta  \Big) =
    O(n^{-40s}).
\end{equation}
Therefore
\begin{equation}
   \sin  \Big(  \frac{2\pi m [\theta N]}{P_{\br-\balpha}}  +\eta\Big) -
   \bE_{\btheta} \Big( \sin  \Big(  \frac{2\pi m [\theta N]}{P_{\br-\balpha}} +\eta  \Big)\; |\; \cF_k \Big) =
    O(n^{-40s}).
\end{equation}
\textsf{ Consider}   \eqref{Le9-4}. By \eqref{Le9-0} and \eqref{Le9-10}, we obtain
\begin{multline*}
  \big| {\bf E}_{\btheta} \big(\sin(2\pi  m[\theta N]/P_{\br -\balpha}
\; | \; \cF_{k-1} )\big) \big|    \leq      \min \Big( 2,  \frac{2^{l_{k-1}-1} }{N \llangle m/P_{\br -\balpha} \rrangle}\Big) \\
=\min \Big( 2,  \frac{P_{\br -\balpha}  2^{l_{k-1}-1} }{m N }\Big) \ll n^{-40s}.
\end{multline*}
\textsf{ Consider}   \eqref{Le9-5}. Let  $ |m_1|/P_{\br_1 - \balpha_1}
\neq   |m_2|/P_{\br_2 -\balpha_2} $.
Applying  Corollary 2, we get
\begin{multline} \nonumber 
\llangle m_1/P_{\br_1 - \balpha_1} +(-1)^j   m_2/P_{\br_2 -\balpha_2} \rrangle =
|  m_1/P_{\br_1 - \balpha_1} +(-1)^j  m_2/P_{\br_2 -\balpha_2}| \\
    > \min(  |m_1/P_{\br_1 - \balpha_1}|, |   m_2/P_{\br_2 -\balpha_2}|)  \exp( -C_1 \ln^2 n) ,
\end{multline}
for $50s^2 p_0 \log_2 n \leq \max_j |r_j| \leq n$, $0 < |m_i| \leq n^{10s}$.
Using \eqref{Le2-8}, \eqref{Le9-2} and \eqref{Le9-10},  we have
\begin{multline}
 2 \Big|\bE_{\btheta}\Big(\sin \big(\frac{2\pi  m_1[\theta N]}{P_{\br_1 - \balpha_1}}   \big)
 \sin \big(\frac{2\pi  m_2[\theta N]}{P_{\br_2 - \balpha_2}}   \big)
  \; | \; \cF_{k-1}\Big)    \\
  -   \bE_{\btheta}\Big(\sin \big(\frac{2\pi  m_1[\theta N]}{P_{\br_1 - \balpha_1}}   \big)
 \sin \big(\frac{2\pi  m_2[\theta N]}{P_{\br_2 - \balpha_2}}   \big)\Big) \Big| \\
   \leq \sum_{j=1}^2 \Big|\bE_{\btheta}\Big(\sin \big(\frac{2\pi  m_1[\theta N]}{P_{\br_1 - \balpha_1}}
    + (-1)^j  \frac{m_2[\theta N]}{P_{\br_2 - \balpha_2}}  +\pi/2 \big)
  \; | \; \cF_{k-1}\Big) \\
- \bE_{\btheta}\Big(\sin \big(\frac{2\pi  m_1[\theta N]}{P_{\br_1 - \balpha_1}}
    + (-1)^j \frac{ m_2[\theta N]}{P_{\br_2 - \balpha_2}}  +\pi/2 \big)
  \Big)  \Big| \\
 \leq \sum_{j=1}^2 \min \Big( 2,  \frac{2^{l_{k-1}} }{N \llangle m_1/P_{\br_1 - \balpha_1} +(-1)^i   m_2/P_{\br_2 -\balpha_2} \rrangle }\Big) \\
 \leq  2^{l_{k-1}}\max_j P_{\br_j - \balpha_j} N^{-1}  \exp( C_1 \ln^2 n) \ll   2^{-\log_2^{10s} n }
  \exp( C_1 \ln^2 n) \ll n^{-40s}. \label{Le9-13}
\end{multline}
Now let  $ |m_1|/P_{\br_1 - \balpha_1}
=   |m_2|/P_{\br_2 -\balpha_2} $.
From \eqref{Le2-8}, \eqref{Le9-3} and \eqref{Le9-10},  we obtain
\begin{multline}
 2 |\bE_{\btheta}\Big(\sin^2 \big( \frac{2\pi  m_1[\theta N]}{P_{\br_1 - \balpha_1}}  \big)
 \; | \; \cF_{k-1}\Big) -  \bE_{\btheta}\Big(\sin^2 \big(\frac{2\pi  m_1[\theta N]}{P_{\br_1 - \balpha_1}}  \big)\Big)     | \\
 = |\bE_{\btheta}(\sin \big( \frac{2\pi  m_1[\theta N]}{P_{\br_1 - \balpha_1}} +\pi/2 \big)
 \; | \; \cF_{k-1}) -  \bE_{\btheta}(\sin \big(\frac{2\pi  m_1[\theta N]}{P_{\br_1 - \balpha_1}} +\pi/2 \big))
      |  \ll n^{-40s}.  \label{Le9-14}
%
\end{multline}
By \eqref{Le9-13} and \eqref{Le9-14}, \eqref{Le9-5} is  proved and Lemma 16 follows. \qed \\ \\

{\bf Bibliography.} \\


[ADDS] Avila, A.,  Dolgopyat, D., Duryev, E.,  Sarig, O.,  The visits to zero of a random walk driven by an irrational
rotation, Israel Journal of Mathematics, 207(2):653-717, 2015.

[Be1] Beck, J., Randomness of the square root of $2$ and the giant leap, Part 1,
Period. Math. Hungar., 60(2):137-242, 2010.

[Be2]
Beck, J., Randomness of the square root of $2$ and the giant leap, Part 2,
Period. Math. Hungar., 62(2):127-246, 2011.

[BeCh]
 Beck, J., Chen, W.~W.~L.,
 Irregularities of Distribution,
 Cambridge Univ. Press,  Cambridge, 1987.

[BrUl] Bromberg, M. and Ulcigrai, U., A temporal central limit theorem
for real-valued cocycles over rotations,
Ann. Inst. Henri Poincar\'e Probab. Stat. 54 (2018), no. 4, 2304-2334.

[FKP]  Faure, H., Kritzer, P., Pillichshammer F.,
From van der Corput to modern constructions of sequences for quasi-Monte Carlo rules,
   Indagationes Mathematicae,  26 (2015),  760-822.

[DS] Dolgopyat, D., Sarig, O., Temporal distributional limit theorems for dynamical systems,
 J. Stat. Phys. 166 (2017), no. 3-4, 680–713.


[Du] Durrett, R.,  Probability: Theory and Examples, Fourth edition,  Cambridge, 2010.

[Ha] Hall, P.,  Heyde, C.C.,
Martingale Limit Theory and Its Application,  New York, 1980.

[Le1]  Levin, M.B., On the upper bound of the $L_p$ discrepancy  of Halton's sequence and the
 Central Limit Theorem for Hammersley's net, ArXiv:1806.11498.

[Le2]
Levin, M.B., On the Gaussian limiting distribution of lattice points in a parallelepiped,
Unif. Distrib. Theory, 11 (2016), no. 2, 45-89.

[LM]
 Levin, M.B., Merzbach, E., Central limit theorems for the ergodic adding machine,
 Israel J. Math., 134 (2003), 61-92.

[Ni]  Niederreiter, H., Random Number Generation and Quasi-Monte Carlo Methods,  SIAM, 1992.

[Wi] Wills, J. M., Zur Gitterpunktanzahl konvexer Mengen,
Elem. Math.,  28 (1973), 57-63.
\\ \\

{\bf Address}: Department of Mathematics,
Bar-Ilan University, Ramat-Gan, 5290002, Israel \\
{\bf E-mail}: mlevin@math.biu.ac.il\\

\end{document}